\documentclass{icm2010}
\usepackage{amscd,amssymb}
\usepackage{graphicx}
\title[Dynamics on geometrically finite hyperbolic manifolds]{Dynamics on geometrically finite hyperbolic manifolds
with applications to Apollonian circle packings and beyond}
\author[Hee Oh]{Hee Oh
}
\contact[heeoh@math.brown.edu]{Mathematics department, Brown university, Providence, RI
and Korea Institute for Advanced Study, Seoul, Korea}

\newtheorem{Thm}[equation]{Theorem}

\newtheorem{Rmk}[equation]{Remark}
\newtheorem{Cor}[equation]{Corollary}

\newtheorem{Def}[equation]{Definition}
\numberwithin{equation}{section}

\newcommand{\z}{\mathbb{Z}}

\renewcommand{\c}{\mathbb{C}}
\newcommand{\br}{\mathbb{R}}

\newcommand{\ba}{\backslash}

\newcommand{\G}{\Gamma}

\renewcommand{\P}{\mathcal P}

\newcommand{\la}{\langle}
\newcommand{\ra}{\rangle}

\newcommand{\SL}{\operatorname{SL}}

\newcommand{\bp}{\begin{pmatrix}}
\newcommand{\ep}{\end{pmatrix}}

\newcommand{\SO}{\operatorname{SO}}

\newcommand{\PSL}{\op{PSL}}
\newcommand{\BMS}{\op{BMS}}
\newcommand{\BR}{\op{BR}}\newcommand{\sk}{\op{sk}}
\newcommand{\bi}{\begin{itemize}}
\newcommand{\be}{\begin{enumerate}}

\newcommand{\ee}{\end{enumerate}}

\newcommand{\op}{\operatorname}

\newcommand{\PS}{\operatorname{PS}}
\newcommand{\Leb}{\operatorname{Leb}}
\newcommand{\bH}{\mathbb H}\newcommand{\LG}{\Lambda_\G}
\newcommand{\T}{\op{T}}
\renewcommand{\PS}{\operatorname{PS}}
\renewcommand{\BR}{\operatorname{BR}}
\renewcommand{\Leb}{\operatorname{Leb}}
\renewcommand{\BMS}{\operatorname{BMS}}
\renewcommand{\LG}{\Lambda(\G)}

\begin{document}

\begin{abstract} We present recent results on counting and distribution of circles
in a given circle packing invariant under a geometrically finite Kleinian group
and discuss how the dynamics of flows on  geometrically finite hyperbolic
$3$ manifolds are related. Our results apply to Apollonian circle packings, Sierpinski curves, Schottky dances, etc.
 \end{abstract}

%\thanks{partially supported by NSF
%   grant 0629322}

% Your AMS 200 Classification should come here %

% Any keywords %

%\maketitle
%\tableofcontents

\section{Introduction}
Let $G$ be a connected semisimple Lie group and $\G <G$ a discrete subgroup with finite co-volume.
Dynamics of flows on the
 homogeneous space $\G\ba G$ have been studied intensively over the last several decades
and brought many surprising applications in various fields notably including
analytic number theory, arithmetic geometry and Riemmanian geometry (see
\cite{MargulisICM1991}, \cite{RatnerICM1994}, \cite{DaniICM1994}, \cite{EskinICM1998},
\cite{KleinbockMargulis1998}, \cite{Yukie1997}, \cite{LubotzkyZimmer2001}, \cite{UllmoICM2002},
\cite{ElkiesMcMullen2004}, \cite{EskinOh2006}, \cite{EskinOh2006d}, \cite{VatsalICM2006}, \cite{LindenstraussICM2006}, \cite{MichelICM2006}, \cite{KlinglerYafaev}, \cite{UllmoYafaev}, \cite{GorodnikOh2007}, \cite{GorodnikShahOh}, \cite{GorodnikOh2008}, \cite{ShahICM2010}, etc.)
The assumption that the volume of $\G\ba G$ is finite is crucial in
most developments in the ergodic theory for flows on $\G\ba G$,
as many basic ergodic theorems fail in the setting of an infinite measure space.
It is unclear what kind of measure theoretic and topological rigidity for
flows on $\G\ba G$ can be expected for a general discrete subgroup $\G$.

In this article we consider the situation when $G$ is
the isometry group of the real hyperbolic space $\bH^n$, $n\ge 2$,
 and $\G<G$ is a geometrically finite discrete  subgroup.
%However when $G$ is the isometry group of the real hyperbolic space $\bH^n$
% and $\G<G$ is the so-called geometrically finite (torsion-free) group,
In such cases we have a rich theory of
 the Patterson-Sullivan density and
  the structure of a fundamental domain for $\G$ in $\bH^n$ is well understood.
Using these we obtain certain equidistribution results for
specific flows on the unit tangent bundle $\T^1(\G\ba \bH^n)$
and apply them to prove results on counting and equidistribution for
circles in a given circle packing of the plane (and also of the sphere) invariant under geometrically finite
groups.

%In this article we discuss some of natural counting and equidistribution problems for points,
%circles and spheres in Euclidean spaces
% which are related to the dynamics of flows on infinite volume geometrically finite
% hyperbolic manifold.
 There are numerous natural questions which arise
 from the analogy with the finite volume cases and most of them are unsolved. We address some of them
 in the last section. We remark that an article by Sarig \cite{Sarig2010} discusses related issues
 but for geometrically {\it{infinite}} surfaces.

 \bigskip
\noindent{\bf Acknowledgement:}
I would like to thank Peter Sarnak for introducing  Apollonian circle packings to me
and for the encouragement to work on this project. I am grateful to Curt McMullen for showing me
the picture of Sierpinski curve which led me to think about more general circle packings
beyond Apollonian ones, as well as for many valuable discussions. I thank my collaborators
Nimish Shah and Alex Kontorovich for the joint work. I also thank Marc Burger and Gregory Margulis
 for carefully reading an earlier draft and making many helpful comments.
Finally I thank my family for their love and support always.

 \section{Preliminaries}
We review some of basic definitions as well as set up notations.
 Let $G$ be the identity component of the isometry group
of the real hyperbolic space $\bH^n$, $n\ge 2$.
Let $\G<G$ be a torsion-free discrete subgroup.
We denote by $\partial_\infty(\bH^n)$ the geometric boundary of $\bH^n$. The limit set
  $\Lambda(\G)$ of $\G$ is defined to be the
 set of accumulation points of an orbit of $\G$ in $\bH^n\cup \partial_\infty(\bH^n)$.
As $\G$ acts on $\bH^n$ properly discontinuously, $\Lambda(\G)$ lies in $\partial_\infty(\bH^n)$.
 Its complement $\Omega(\G):=\partial_\infty(\bH^n) -\Lambda(\G )$ is called the domain of discontinuity for $\G$.

An element $g\in G$ is called parabolic if it fixes a unique point in  $\partial_\infty(\bH^n)$ and
loxodromic if it fixes two points in  $\partial_\infty(\bH^n)$.
A limit point $\xi\in \Lambda(\G)$ is called a parabolic fixed point
if it is fixed by a parabolic element of $\G$ and called a radial limit point
 (or a conical limit point or  a point of
approximation) if for some geodesic ray
$\beta$ tending to $\xi$ and some point
$x\in \bH^n$, there is a sequence $\gamma_i\in \G$ with
$\gamma_i x\to \xi$ and $d(\gamma_i x, \beta)$ is bounded,
 where $d$ denotes the hyperbolic
distance.
A parabolic fixed point $\xi$ is called bounded if $\op{Stab}_\G(\xi)\ba (\Lambda(\G)-\{\xi\})$
is compact.

The convex core $C_\G$ of $\G$ is defined to be the minimal convex
set in $\bH^n$ mod $\G$ which contains all geodesics connecting
any two points in $\Lambda(\G)$.
A discrete subgroup $\G$ is called {\it  geometrically finite}
 if the unit neighborhood of its convex core has finite volume and called {\it convex co-compact}
 if its convex core is compact.
 It is clear that
 a (resp. co-compact) lattice in $G$ is geometrically finite (resp. convex co-compact).
Bowditch showed \cite{Bowditch1993} that $\G$ is geometrically finite if and only if
$\Lambda(\G)$ consists entirely of radial limit points and bounded parabolic fixed points.
It is further equivalent to saying that $\G$ is finitely generated for $n=2$, and that
$\G$ admits a finite sided fundamental domain in $\bH^3$ for $n= 3$.
We refer to \cite{Bowditch1993} for other equivalent definitions.

$\G$ is called {\it elementary} if $\Lambda(\G)$ consists of at most two points, or
 equivalently, $\G$ has an abelian subgroup of finite index.

We denote by $0\le \delta_\G\le n-1 $ the critical exponent of $\G$, that is,
the abscissa of convergence
of the Poincare series of $\G$: $$\mathcal P_\G(s)
:=\sum_{\gamma\in\G} e^{-s d(o, \gamma o)}$$
where $o\in \bH^n$. For a non-elementary group  $\G$, $\delta_\G$ is positive and Sullivan
\cite{Sullivan1979} showed that for $\G$ geometrically finite, $\delta_\G$ is equal to
  the Hausdorff dimension of the limit set $\Lambda(\G)$.

For $\xi\in \partial_\infty(\bH^n)$ and $y_1, y_2\in \bH^n$, the Busemann function $\beta_\xi(y_1, y_2)$
 measures a signed distance between horospheres
passing through  $y_1$ and $y_2$ based at $\xi$:
$$\beta_\xi(y_1,y_2)=\lim_{t\to\infty} d(y_1, \xi_t)-d(y_2,\xi_t)$$ where
$\xi_t$ is a geodesic ray toward  $\xi$.

 For a vector $u$ in the unit tangent bundle $\T^1(\bH^n)$, we define $u^{\pm}\in \partial_\infty(\bH^n)$ to be the
 two end points of the geodesic determined by $u$:
 $$u^+:=\lim_{t\to  \infty} g^t(u)\quad\text{and }\quad  u^-
 :=\lim_{t\to -\infty} g^t(u) $$
where $\{g^t\}$ denotes the geodesic flow.

We denote by $\pi: \T^1(\bH^n)\to \bH^n$ the canonical projection.
Fixing a base point $o\in \bH^n$, the map $$u\mapsto (u^+, u^-,
\beta_{u^-}(\pi(u), o))$$ yields a homeomorphism between
$\T^1(\bH^n)$ and $(\partial_\infty(\bH^n)\times \partial_\infty(\bH^n) - \{(\xi,\xi):\xi\in
\partial_\infty(\bH^n)\}) \times \br$.

 \bigskip

Throughout the paper we assume that $\G$ is non-elementary.

\bigskip

\noindent{\bf {Patterson-Sullivan density}:}
Generalizing the work of Patterson \cite{Patterson1976} for $n=2$,
Sullivan \cite{Sullivan1979} constructed a $\G$-invariant conformal density $\{\nu_x: x\in \bH^n\}$
of dimension $\delta_\G$ on $\Lambda(\G)$.
 That is, each
$\nu_x$ is a finite Borel measure on $\partial_\infty(\bH^n)$ supported on $\Lambda(\G)$
satisfying that for any $x,y\in \bH^n$, $\xi\in \partial_\infty(\bH^n)$ and $\gamma\in
\G$,
$$\gamma_*\nu_x=\nu_{\gamma x}\quad\text{ and}\quad
\frac{d\nu_y}{d\nu_x}(\xi)=e^{-\delta_\G \beta_{\xi} (y,x)}, $$
where $\gamma_*\nu_x(R)=\nu_x(\gamma^{-1}(R))$.

For $\G$ geometrically finite, such conformal density $\{\nu_x\}$ exists uniquely up to homothety.
In fact, fixing $o\in
\bH^n$,  $\{\nu_x\}$ is a constant multiple
of the following family $\{\nu_{x,o}\}$ where $\nu_{x,o}$ is the weak-limit as $s\to \delta_\G^+$
of the family of measures
$$\nu_{x,o}(s):=\frac{1}{\sum_{\gamma\in \G} e^{-sd(o, \gamma o)}}
\sum_{\gamma\in\G} e^{-sd(x, \gamma o)} \delta_{\gamma o} $$
where $\delta_{\gamma o}$ denotes the Dirac measure at $\gamma o$.

Consider the Laplacian $\Delta$ on $\bH^n$. In the upper half-space coordinates
$\bH^n=\{(x_1, \cdots, x_{n-1}, y):y>0\}$ with the metric
$\frac{\sqrt{dx_1^2+\cdots +dx_{n-1}^2+dy^2}}{y}$,
it is given as $$\Delta=-y^2\left(\frac{\partial^2}{\partial x_1^2}+\cdots
 +\frac{\partial^2}{\partial x_{n-1}^2} +\frac{\partial^2}{\partial y^2}\right) +(n-2)y\frac{\partial}{\partial y} $$
 (strictly speaking, this is the negative of the usual hyperbolic Laplacian).
 Sullivan  \cite{Sullivan1979} showed that
 $$\phi_\G (x):=|\nu_x|$$ is an eigenfunction
for $\Delta$ with eigenvalue $\delta_\G(n-1-\delta_\G)$.
From the $\G$-invariance of the Patterson-Sullivan density $\{\nu_x\}$, $\phi_\G$ is a function on $\G\ba \bH^n$.
Sullivan further showed that if $\G$ geometrically finite and $\delta_\G >(n-1)/2$,
$\phi_\G$ belongs to $L^2(\G\ba \bH^n)$ and is a unique (up to a constant multiple)
positive eigenfunction with the smallest eigenvalue $\delta_\G(n-1-\delta_\G)$
(cf. \cite{Sullivan1987}). Combined with a result of Yau \cite{Yau1975}, it follows that
  $\delta_\G=n-1$ if and only if $\G$ is a lattice in $G$.

 \bigskip

%\begin{Def}[Bowen-Margulis-Sullivan measure]\label{defbms}
%{\rm
\noindent{\bf {Bowen-Margulis-Sullivan measure}:}
Fixing the Patterson-Sullivan density $\{\nu_x\}$,
the Bowen-Margulis-Sullivan measure $m^{\BMS}_\G$
(\cite{Bowen1971}, \cite{Margulisthesis}, \cite{Sullivan1984}) is the induced measure on $\T^1(\G\ba \bH^n)$ of the following
$\G$-invariant measure on $\T^1(\bH^n)$:
$$d  \tilde m^{\BMS}(u)=e^{\delta_\G \beta_{u^+}(x, \pi(u))}\;
 e^{\delta_\G \beta_{u^-}(x,\pi(u)) }\;
d\nu_{x}(u^+) d\nu_{x}(u^-) dt $$
 where $x\in \bH^n$.
   %via the homeomorphism between $\T^1(\bH^n)$
%and $(\partial_\infty(\bH^n)\times \partial_\infty(\bH^n) - \{(\xi,\xi):\xi\in
%\partial_\infty(\bH^n)\}) \times \br$ given by
%$u\mapsto (u^+, u^-,
%t:=\beta_{u^-}(\pi(u), x))$, where
% $\pi: \T^1(\bH^n)\to \bH^n$ denotes the canonical projection.
%and $\pi: \T^1(\bH^n) \to \bH^n$ denotes the canonical projection.
%

It follows from the conformality of $\{\nu_x\}$ that
this definition is independent of the choice of $x$.
% By the uniqueness of the Patterson-Sullivan density, $m^{\BMS}_\G$ is uniquely determined
 % up to a constant multiple as well.
 The measure $m^{\BMS}_\G$ is invariant under the geodesic flow
and is supported on
the non-wandering set
$ \{u\in \T^1(\G\ba \bH^n): u^{\pm}\in \Lambda(\G)\}$ of the geodesic flow.
 Sullivan showed that for $\G$ geometrically finite, the total mass
$|m^{\BMS}_\G|$ is finite and the geodesic flow is ergodic with respect to $m^{\BMS}_\G$
\cite{Sullivan1984}.
This is a very important point for the ergodic theory on geometrically finite hyperbolic manifolds,
since despite of  the fact that the Liouville measure is infinite, we do have
  a finite measure on $\T^1(\Gamma\ba \bH^n)$ which is invariant and ergodic for the geodesic flow.
Rudolph \cite{Rudolph1982} showed that the geodesic flow is even mixing with respect to $m^{\BMS}_\G$.

\section{Counting and distribution of circles in the plane}
A circle packing in the plane $\c$ is simply a union of circles.
As circles may intersect with each other beyond tangency points, our definition
of a circle packing is more general than what is usually thought of.
For a given circle packing $\P$ in the plane, we discuss
questions on counting and distribution of small circles in $\P$.
A natural size of a circle is measured by its radius. We will use the curvature (=the reciprocal of
the radius) of a circle instead.

\begin{figure}
 \begin{center}
    \includegraphics[width=5cm]{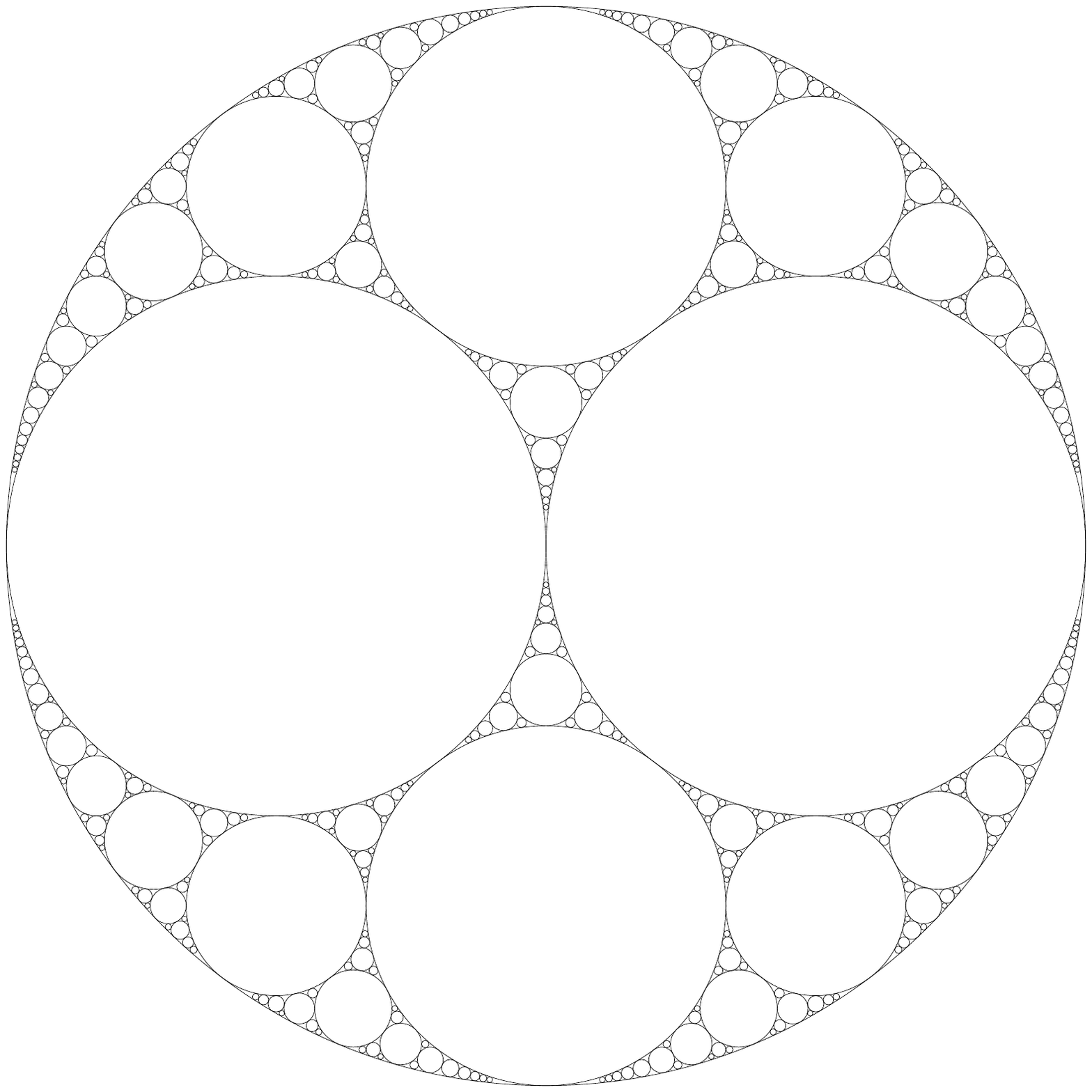}
    \includegraphics[width=5cm]{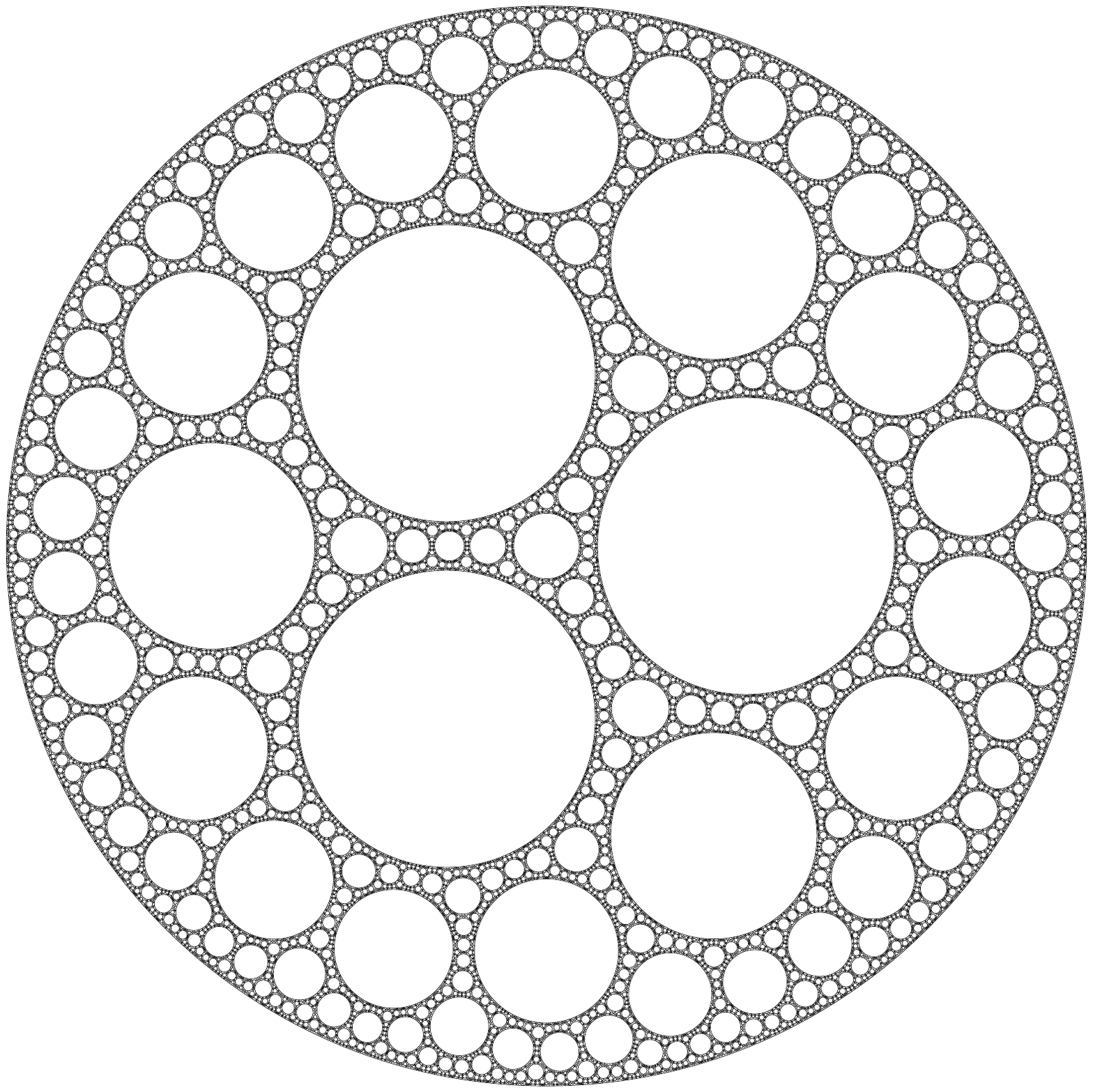}
    \caption{Apollonian circle packing and Sierpinski curve (by C. McMullen)}
    \label{ASplane}
 \end{center}
\end{figure}

We suppose that $\P$ is infinite and that
 $\mathcal P$ is locally finite in the sense that
for any $T>0$, there
are only finitely many circles of curvature at most $T$ in any fixed bounded region of the plane.
See Fig. \ref{ASplane}, \ref{fractal} and \ref{fpar} for examples of locally finite packings.

For a bounded region $E$ in the plane $\c$,
we consider
the following counting function:
$$N_T(\P, E):=\# \{C\in \P: C \cap
E\ne\emptyset ,\;\;\op{Curv}(C)<T\}  $$
where $\op{Curv}(C)$ denotes the curvature of $C$.
The local finiteness assumption is so that $N_T(\P, E)<\infty $ for any bounded $E$.
We ask if there is an asymptotic for $N_T(\P, E)$ as $T$ tends to infinity
and what the dependence of such an asymptotic on $E$ is, if exists.

Consider the upper half space model
$\bH^3=\{(z,r):z\in \c, r>0\}$ with the hyperbolic metric given by $\frac{\sqrt{|dz|^2+dr^2}}{r}$.
An elementary but helpful observation is that
if we denote by $\hat C\subset \bH^3$ the convex hull
of $C$, that is, the northern hemisphere above $C$,
then
$N_T(\P, E)$ is equal to the number of hemispheres of height at most $T^{-1}$ in $\bH^3$
whose boundaries lie in $\P$ and intersect $E$, as the radius of a circle is same
as the height of the corresponding hemisphere.

Let $\G<\PSL_2(\c)$ be a geometrically finite discrete subgroup
and fix a $\G$-invariant Patterson-Sullivan density $\{\nu_x:x\in \bH^3 \}$.

In order to present our theorem on the asymptotic of $N_T(\P, E)$
for $\P$ invariant under $\G$, we introduce two new invariants associated to $\G$ and $\P$.
The first one is a Borel measure on $\c$ which depends only on $\G$.

\begin{Def}{\rm  Define
a Borel measure $\omega_\G$ on $\c$: for $\psi\in C_c(\c)$,
$$\omega_\G(\psi)=\int_{z\in \c} \psi(z) e^{\delta_\G \beta_z(x,z+j) }\; d\nu_{x}(z) $$
where $j=(0,1)\in \bH^3$ and $x\in \bH^3$.
By the conformal property of $\{\nu_{x}\}$, this definition is independent of the choice
of $x\in \bH^3$.} \end{Def}

 Note that $\omega_\G$ is supported on $\Lambda(\G)\cap \c$
 and in particular that $\omega_\G(E)>0$ if the interior of $E$
 intersects $\Lambda(\G)\cap \c$ non-trivially.
 We compute:
$$d\omega_{\G} =(|z|^2+1)^{\delta_\G} d\nu_{j} .$$

The second one is a number in $[0,\infty]$ measuring certain size of $\P$:
\begin{Def}[The $\G$-skinning size of $\P$]\label{sk}
 {\rm For a circle packing $\P$ invariant under $\G$, we define:
$$\op{sk}_\G(\P):=\sum_{i\in I} \int_{s\in \op{Stab}_{\G} (C_i^\dagger)\ba C_i^\dagger}  e^{\delta_\G
\beta_{s^+}(x,\pi(s))}d\nu_{x}(s^+)$$
where $x\in \bH^3$, $\{C_i:i\in I\}$ is a set of representatives of $\G$-orbits in $\P$ and
$C_i^\dagger\subset \T^1(\bH^3)$ is the set of unit normal vectors to the convex hull $\hat C_i$ of $C_i$.
 Again by the conformal property of $\{\nu_x\}$, the definition of
 $\sk_\G(\P)$ is independent of the choice of $x$ and the choice of representatives $\{C_i\}$. }\end{Def}

We remark that the value of $ \op{sk}_\G(\P)$ can be zero
or infinite in general and we do not assume any condition
on $\op{Stab}_{\G} (C_i^\dagger)$'s (they may even be trivial).
%For instance, if $\Lambda(\G)$ is precisely a circle $C$,
%then $\{s\in C^\dag: s^{\pm}\in \Lambda(\G)\}=\emptyset$ and hence
% $\op{sk}_\G(\G(C))=0$.
By the interior of a circle $C$, we mean the open disk which is enclosed by $C$.
We then have the following:
\begin{Thm}[\cite{OhShahcircle}]\label{m1}
Let $\G$ be a non-elementary geometrically finite discrete subgroup of $\PSL_2(\c)$ and let
 $\mathcal P=\cup_{i\in I}\Gamma(C_i)$ be an infinite,
  locally finite, and  $\G$-invariant circle packing with finitely many $\G$-orbits.

Suppose one of the following conditions hold:
\begin{enumerate}
\item $\G$ is convex co-compact;
\item  all circles in $\P$ are mutually disjoint;
\item $\cup_{i\in I}C_i^\circ \subset \Omega(\G)$
where $C_i^\circ$ denotes the interior of $C_i$.
\end{enumerate}For any bounded region $E$ of $\c$
whose boundary is of zero Patterson-Sullivan measure,
%a piecewise algebraic curve,
 we have
 $$N_T(\P, E)
\sim
 \frac{ \sk_\G(\P) }{\delta_\G \cdot |m^{\BMS}_\G|} \cdot
\omega_\G(E)  \cdot T^{\delta_\G} \quad \text{ as $T\to \infty$}$$
 and $0<\sk_\G(\P)<\infty$.
 \end{Thm}
 \begin{Rmk}{\rm
\begin{enumerate}
\item If $\G$ is Zariski dense in $\PSL_2(\c)$,
considered as a real algebraic group, any
real algebraic curve has
 zero Patterson-Sullivan measure \cite[Cor. 1.4]{FlaminioSpatzier}.
Hence the above theorem applies to any Borel subset $E$ whose
boundary is a countable union of real algebraic curves.
 \item  We call
the complement in $\hat \c$ of the set $\cup_{i\in I} \G (C_i^\circ)$ the residual set of $\P$.
The condition (3) above is then equivalent to saying that $\Lambda(\G)$ is contained in the residual set of $\P$.

\item
If we denote by $H^-_\infty(j)$ the contracting horosphere based at $\infty$ in $\T^1(\bH^3)$
which
consists of all upward normal unit vectors on $\c +j=\{(z,1):z\in \c\}$,
we can alternative write
the measure $\omega_\G$ as follows:
$$\omega_\G(\psi)=\int_{u\in H^-_\infty(j)} \psi(u^-) e^{\delta_\G
\beta_{u^-}(x, \pi(u)) }\; d\nu_{x}(u^-) $$
and recognize that $\omega_\G$ is the projection
of the conditional of the Bowen-Margulis-Sullivan measure
$\tilde m^{\BMS}$ on the horosphere $H^-_\infty(j)$ to $\c$ via the map
$u\mapsto u^-$. It is worthwhile to note that
the hyperbolic metric on $\c+j$ is precisely the Euclidean metric.

 \item
Suppose that circles in $\P$ are disjoint possibly except for tangency points
and that $\Lambda(\G)$ is equal to the residual set of $\P$.
If $\infty$ is either in $\Omega(\G)$ (that is, $\P$ is bounded)
or a parabolic fixed point
for $\G$, then
 $\delta_\G$ is equal to the circle packing exponent $e_{\mathcal P}$ given by
 $$e_{\mathcal P}=\inf\{s: \sum_{C\in \P} r(C)^s<\infty\}=\sup\{s:\sum_{C\in \P}r(C)^s=\infty\} $$
 where $r(C)$ denotes the radius of $C$ \cite{Parker1995}. This
 extends the earlier work of Boyd \cite{Boyd1973} on bounded
 Apollonian circle packings.
\end{enumerate}}\end{Rmk}

\begin{figure}
\begin{center}
    \includegraphics[width=4in]{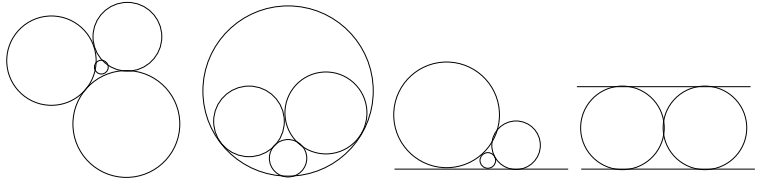}
 \end{center}
 \caption{Possible configurations of four mutually tangent circles}
  \label{Fourcircles}
  \end{figure}

%We now discuss when we can count circles in
 %a general circle packing $\P$, simply a union of circles in the plane.

%We are interested in counting circles in $\mathcal P$ of curvature (=reciprocal of radius)
%at most $T$, as $T\to \infty$. For this,

%For a bounded subset $E$ of $\c$ and $T>0$, we set
%$$N_T(\P, E):=\# \{C\in \P: C\cap
%E\ne\emptyset ,\;\;\op{Curv}(C)<T\} . $$
%In other words, $\P$ is locally finite if and only if
%$N_T(\P, E)<\infty$ for any $T>0$.

We discuss some concrete circle packings to which our theorem applies.
\subsection{Apollonian circle packings in the plane}
Apollonian circle packings are one of the most beautiful circle packings whose
construction can be described in a very simple manner based on
an old theorem of Apollonius (262-190 BC).
It says that given three mutually tangent circles in the plane,
there are exactly two circles which are tangent to all the three circles.

In order to construct an Apollonian circle packing,
we start with four mutually tangent circles. See Fig. \ref{Fourcircles} for possible configurations.
By Apollonius' theorem, there are precisely four new circles that are
tangent to three of the four circles. Continuing to repeatedly add new circles tangent to
three of the circles from the previous generations, we arrive at an infinite circle packing,
 called an Apollonian circle packing,

See Fig.  \ref{boundedApo} and \ref{fpar} for examples of Apollonian circle  packings
where each circle is labeled by its curvature (that is, the reciprocal of its radius).
There are also Apollonian packings which spread all over the plane as well as
spread all over to the half plane. As circles in these packings
 would become enormously large after a few first generations,
it is harder to draw them on paper.

There are many natural questions about Apollonian circle packings
either from the number theoretic or the geometric point of view and we refer to
the series of papers by Graham, Lagarias, Mallows, Wilks, and Yan especially
\cite{GrahamLagariasMallowsWilksYanI} \cite{GrahamLagariasMallowsWilksYanI-n}, and \cite{ErickssonLagarias2007}
as well as the
 letter of Sarnak to Lagarias \cite{SarnakToLagarias} which inspired the author to work
 on the topic personally.  Also see a more recent article \cite{SarnakMAA}.
\begin{figure}\begin{center}
 \includegraphics [width=2in]{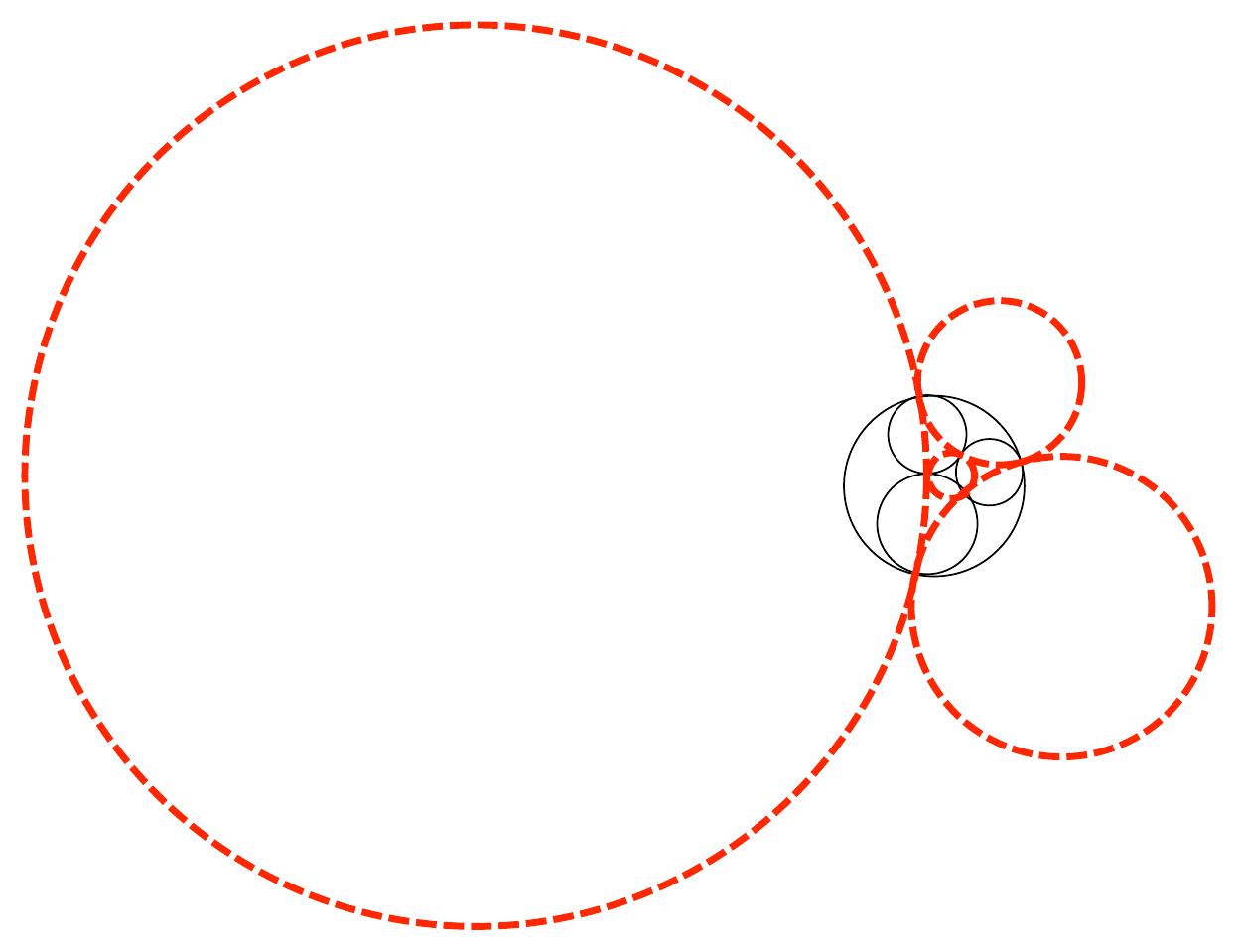}
\caption{Dual circles}
\label{picDual}\end{center}
\end{figure}

To find the symmetry group of a given Apollonian packing $\P$, we
consider the dual circles
%$\tilde C_1, \tilde C_2, \tilde C_3, \tilde C_4$
 to any fixed four mutually
tangent circles
%, say, $C_1, C_2, C_3, C_4$;
(see Fig. \ref{picDual} where the red dotted circles are
the dual circles to the black solid circles).
Inversion with respect to each dual circle
 fixes three circles that the dual circle crosses perpendicularly and interchanges
two circles tangent to those three circles.
Hence the group, say, $\G(\P)$, generated by the four inversions with respect to the dual circles
preserves the packing $\P$ and there are four $\G(\P)$ orbits of circles in $\P$.

As the fundamental domain of $\G(\P)$ in $\bH^3$ can be taken to
be the exterior of the four hemispheres above the dual circles in $\bH^3$,
$\G(\P)$ is geometrically finite.
It is known that the limit set of $\G(\P)$ coincides
 precisely with the residual set of $\P$ and hence
  the critical exponent of $\G(\P)$ is equal to
the Hausdorff dimension of the residual set of $\P$, which
is approximately $$\alpha=1.30568(8)$$ due to C. McMullen \cite{McMullen1998}
 (note that as any two Apollonian packings are equivalent to each other by a Mobius transformation,
$\alpha$ is independent of $\P$). In particular it follows that $\G(\P)$
is Zariski dense in the real algebraic group
$\PSL_2(\c)$ and hence we deduce the following from Theorem \ref{m1} and the remark following it:

\begin{figure}
\begin{center}
\includegraphics [width=4cm]{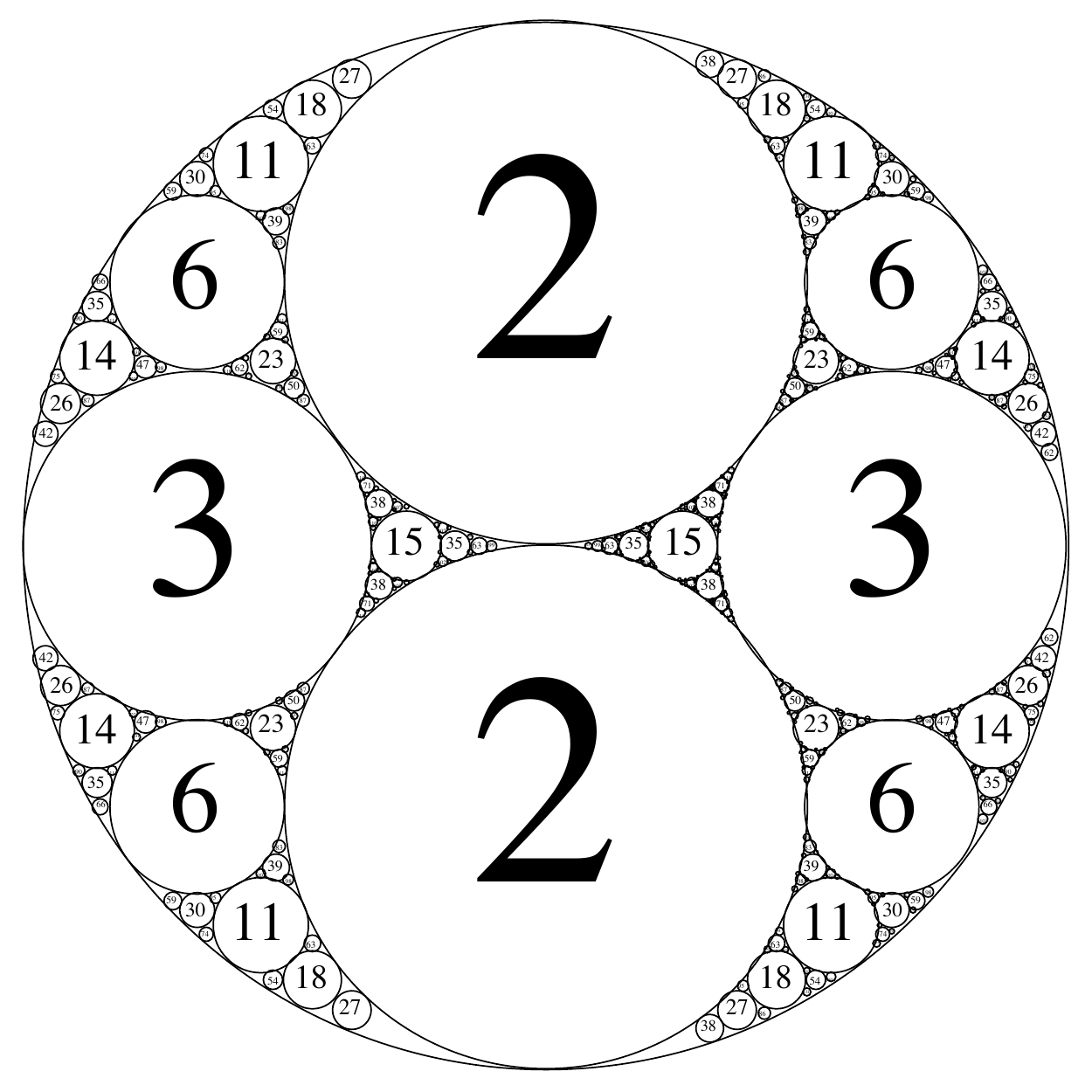}
  \includegraphics[width=4cm]{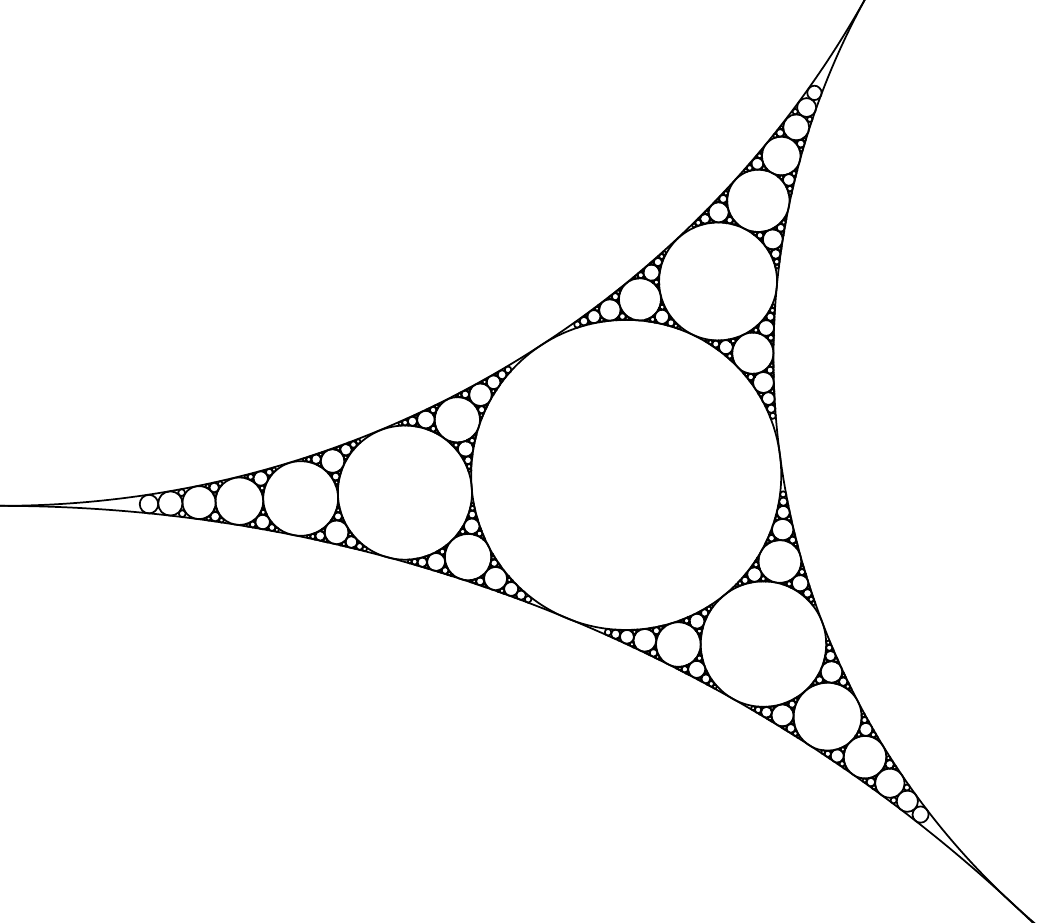}
\end{center}
\caption{A bounded Apollonian circle packing and the Apollonian packing
of a triangular region}
\label{boundedApo}\label{Tri}
\end{figure}

\begin{Cor}[\cite{OhShahcircle}]\label{apo} Let $\P$ be an Apollonian circle packing.
For any bounded region $E$ of $\c$ whose boundary is a countable union of
real algebraic curves,
 we have
 $$N_T(\P, E)
\sim \frac{\sk_{\G_{\P}}(\P) }{\alpha \cdot |m^{\BMS}_{\G_\P}|} \cdot
\omega_{\G_\P}(E)  \cdot T^{\alpha} \quad\text{as $T\to \infty$}$$
where $\G_\P:=\G(\P)\cap \PSL_2(\c)$.
% and $0<\sk_{\G_{\P}}(\P) <\infty$.
% Again by the conformal property of $\{\nu_\xi\}$, the definition of
% $\kappa(\P)$ is independent of the choice of $\xi$ and the choice of representatives $\{C_i\}$.
 \end{Cor}

\begin{Rmk}{\rm\begin{enumerate}
\item In the cases when $\P$ is bounded and $E$ is the largest disk in such $\P$, and when
$\P$ lies between two parallel lines and $E$ is the whole period (see Fig. \ref{fpar}),
the above asymptotic was  previously obtained in \cite{KontorovichOh} with a less explicit
description of the main term.
\item Corollary \ref{apo} applies to any triangular region $\mathcal T$ (see Fig. \ref{Tri}) of an Apollonian circle packing.
\end{enumerate}}
\end{Rmk}

\subsection{More circle packings}
\begin{figure}\begin{center}
 \includegraphics [width=4cm]{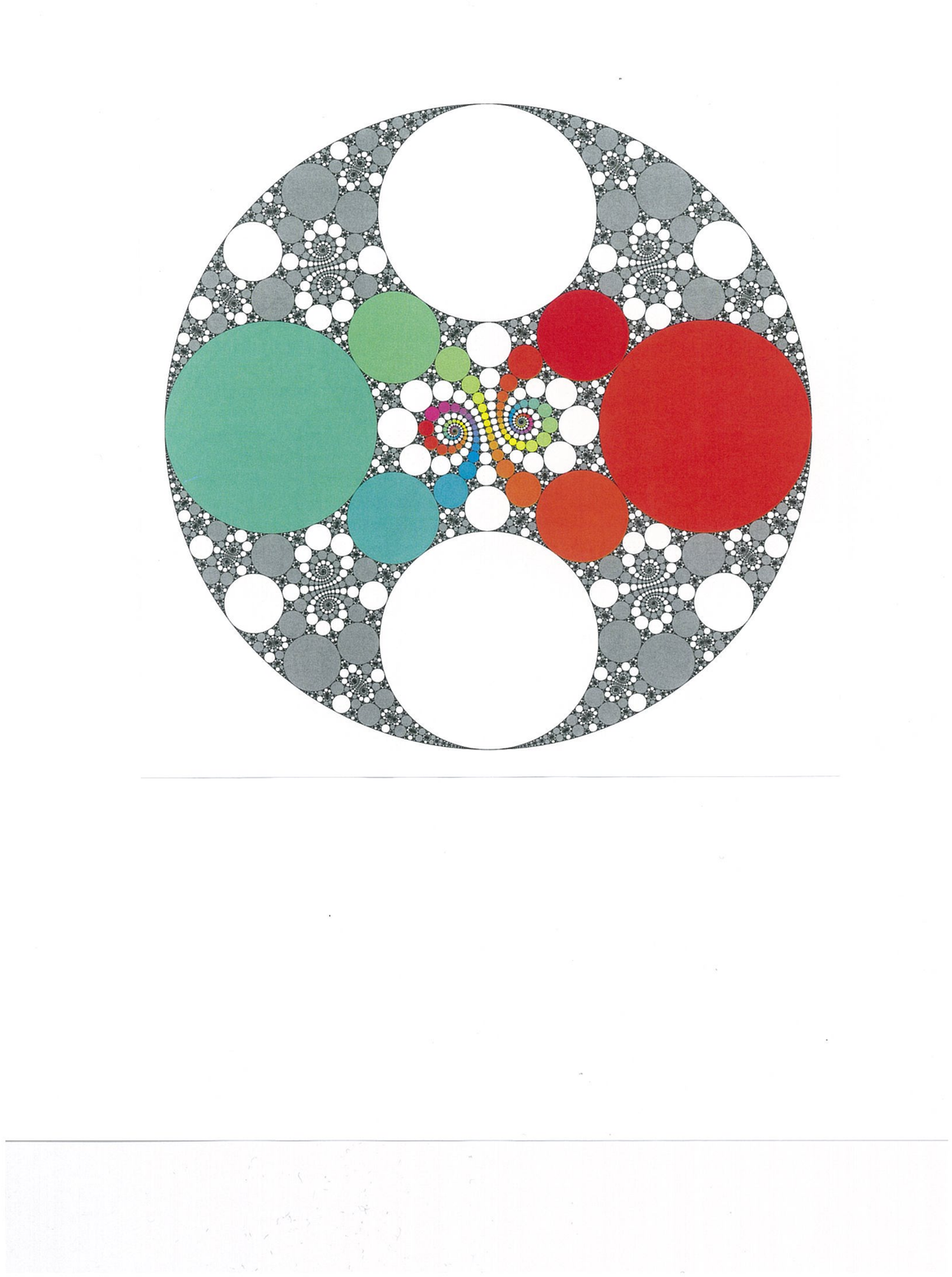}
\includegraphics [width=4cm]{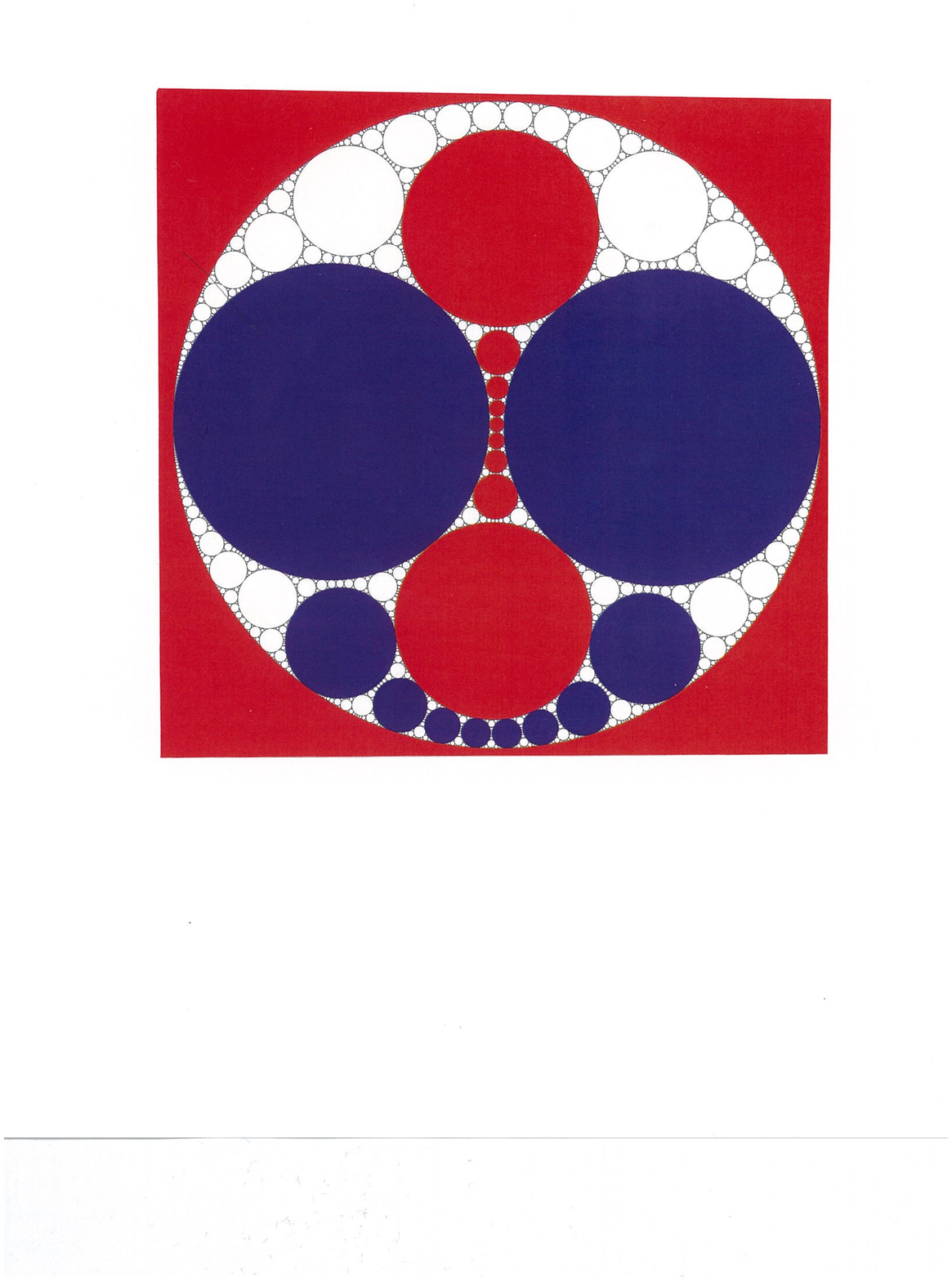}
 \includegraphics [width=4cm]{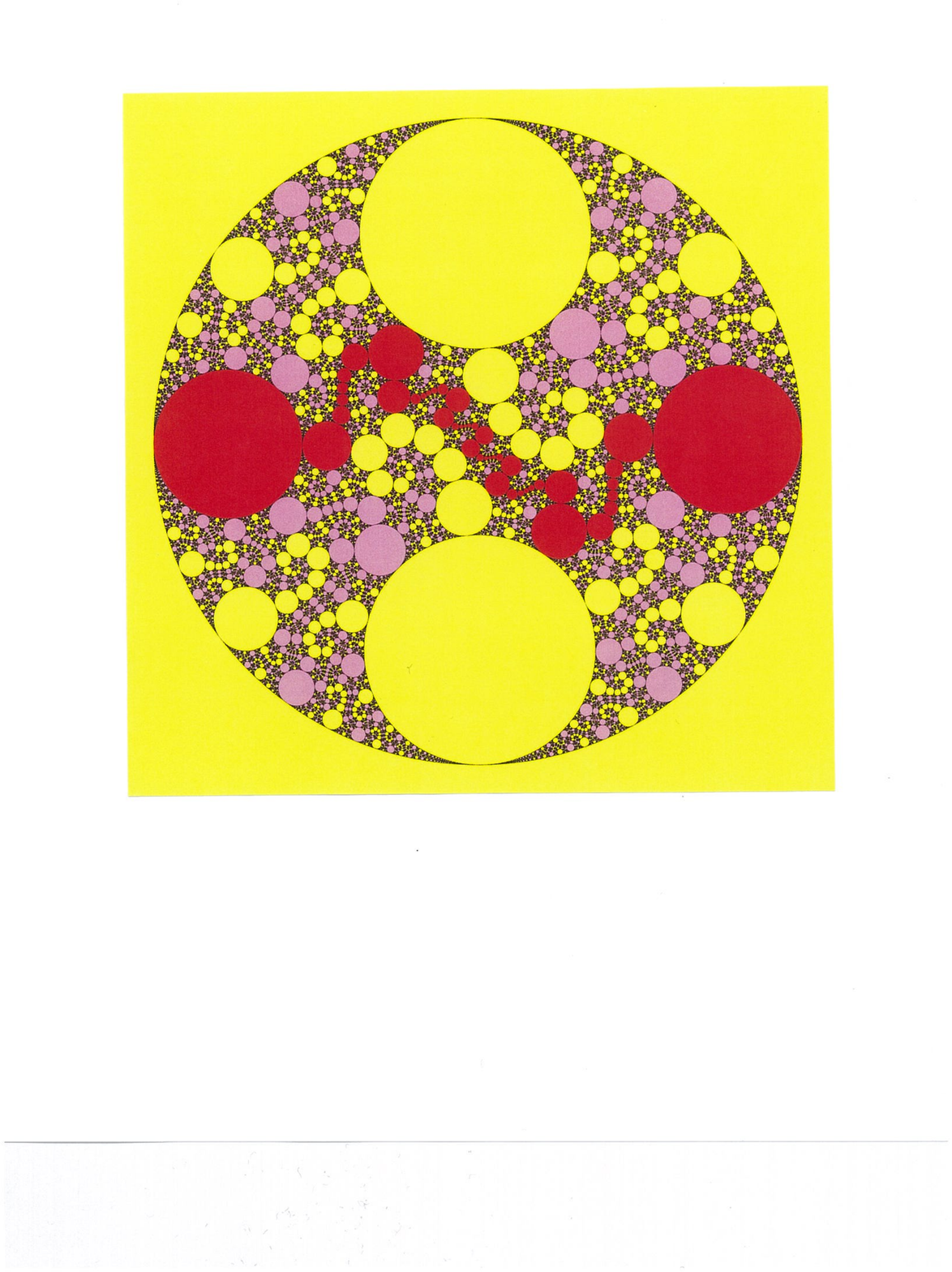}
 \caption{Limit sets of Schottky groups (reproduced with permission from Indra's Pearls, by
D.Mumford, C. Series and D. Wright, copyright Cambridge University Press 2002).} \label{Indra} \end{center}
 \end{figure}
 \subsubsection{Counting circles in the limit set $\Lambda(\G)$}
If $\G\ba \bH^3$ is a hyperbolic $3$ manifold with boundary being totally geodesic, then
$\G$ is automatically geometrically finite \cite{Kojima1992} and $\Omega(\G)$
 is a union of countably many disjoint open disks. Hence Theorem \ref{m1} applies to counting
 these open disks in $\Omega(\G)$ with respect to the curvature, provided there are infinitely many such.
The picture of a Sierpinski curve in Fig. \ref{ASplane} is a special case of this (so
are Apollonian circle packings).
More precisely,  if $\G$ denotes the group generated by reflections in the sides of a unique
 regular tetrahedron whose convex core is bounded by four $\frac{\pi}4$ triangles and
 by four right hexagons, then the residual set of
  a Sierpinski curve in Fig. \ref{ASplane}
 coincides with
  $\Lambda(\G)$ (see \cite {McMullennotes275} for details), and it is known to be homeomorphic to
 the well-known Sierpinski carpet by a theorem of Claytor \cite{Claytor1934}.

% The second picture in Fig. \ref{f3} is taken from McMullen's course notes (see \cite[P.9]{McMullennotes275}).

Three pictures in Fig. \ref{Indra} can be found
 in the beautiful book {\it Indra's pearls} by Mumford, Series and Wright
  \cite{MumfordSeriesWright} and the residual sets are the limit sets
of some (geometrically finite) Schottky groups and hence our theorem applies to counting circles
in those pictures.

\subsubsection{Schottky dance} 
\begin{figure}
 \begin{center}
    \includegraphics[width=5cm]{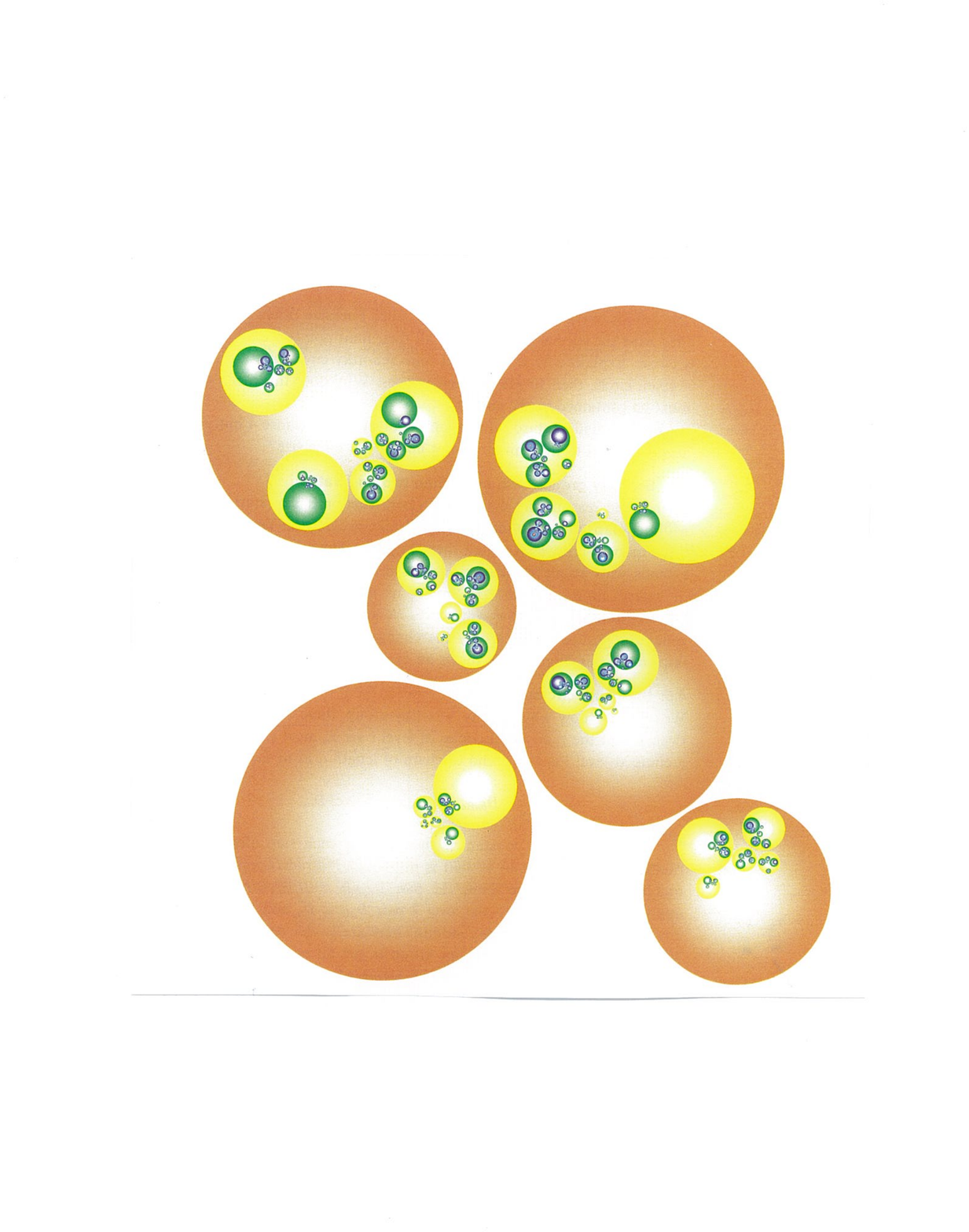}
     \caption{Schottky dance  (reproduced with permission from Indra's Pearls, by
D.Mumford, C. Series and D. Wright, copyright Cambridge University Press 2002)} \label{fractal}
 \end{center}
\end{figure}

Other kinds of examples are obtained by considering the images of
Schottky disks under Schottky groups.
Take $k\ge 1$ pairs of mutually disjoint closed disks
$\{(D_i, D_i'): 1\le i\le k\}$ in $\c$ and choose
 M\"obius transformations $\gamma_i$ which maps $D_i$ and $D_i'$ and sends the interior of
$D_i$ to the exterior of $D_i'$, respectively. The group, say, $\G$,
 generated by $\{\gamma_i:1\le i\le k \}$ is called a Schottky group of genus $k$ (cf. \cite[Sec. 2.7]{MardenOutercircles}).
The $\G$-orbits of the disks nest down onto the limit set $\Lambda(\G)$ which is totally disconnected.
If we denote by $\P$ the union $\cup_{i=1}^k (\G(C_i)\cup \G(C_i'))$
where $C_i$ and $C_i'$ are the boundaries of $D_i$ and $D_i'$ respectively,
$\P$ is locally finite, as the nesting disks will become smaller and smaller (cf. \cite[4.5]{MumfordSeriesWright}).
The common exterior of hemispheres above the initial disks $D_i$ and $D_i'$, $1\le i\le k$,
is a fundamental domain for $\G$ in the upper half-space model $\bH^3$, and hence
$\G$ is geometrically finite.
Since $\P$ consists of disjoint circles,
Theorem \ref{m1} applies to $\P$.
For instance, see Fig. \ref{fractal} (\cite[Fig. 4.11]{MumfordSeriesWright}).
One can find many more explicit circle packings in \cite{MumfordSeriesWright} to which Theorem \ref{m1} applies.

%The following theorem shows that the measure $\omega_{\G}$ describes the asymptotic distribution
%of circles of curvature at most $T$ for $\G$-orbits of circles as $T\to \infty$.

\section{Circle packings on the sphere}
In the unit sphere $\mathbb S^2=\{x^2+y^2+z^2=1\}$
with the Riemannian metric
induced from $\br^3$, the distance between two points is simply the angle between the rays connecting
them to the origin $o=(0,0,0)$.

Let $\P$ be a circle packing on the sphere ${\mathbb S^2}$, i.e., a union of circles.
The spherical curvature of a circle $C$ in $\mathbb S^2$
is given by
$$\op{Curv}_{S}(C)=\cot \theta(C)$$
where $0<\theta(C)\le \pi/2$ is the spherical radius of $C$, that is, the half of the visual angle of $C$
from the origin $o$.
We suppose that $\P$ is infinite and locally finite in the sense that
 there
are only finitely many circles in $\P$ of spherical curvature at most $T$ for any fixed $T>0$.

For a region $E$ of ${\mathbb S^2}$,  we set
$$N_T(\P, E):=\# \{C\in \P: C \cap
E\ne\emptyset ,\;\;\op{Curv}_S(C)<T\} . $$
We consider the Poincare ball model $\mathbb B=\{x_1^2+x_2^2+x_3^2<1\}$
 of the hyperbolic $3$ space with the metric $d$ given by
$\frac{2\sqrt{dx_1^2+dx_2^2+dx_3^2}}{1-(x_1^2+x_2^2+x_3^2)}.$
Note that the geometric boundary of $\mathbb B$ is ${\mathbb S^2}$ and that
for any circle $C$ in $\mathbb S^2$,
we have $$\sin \theta(C)=\frac{1}{\cosh d(\hat C, o)} $$
 where $\hat C\subset \mathbb B$ is the convex hull of $C$.
As both $\sin \theta $ and $\cosh d$ are monotone functions for $0\le \theta\le \pi/2$
 and $d \ge 0$ respectively, understanding $N_T(\P, E)$ is equivalent to investigating
  the number of
 Euclidean hemispheres on $\mathbb B$ meeting the ball of hyperbolic radius $T$ based at $o$
   whose boundaries are in $\P$ and intersect $E$.

\begin{figure}
 \begin{center}
 \includegraphics[width=5cm]{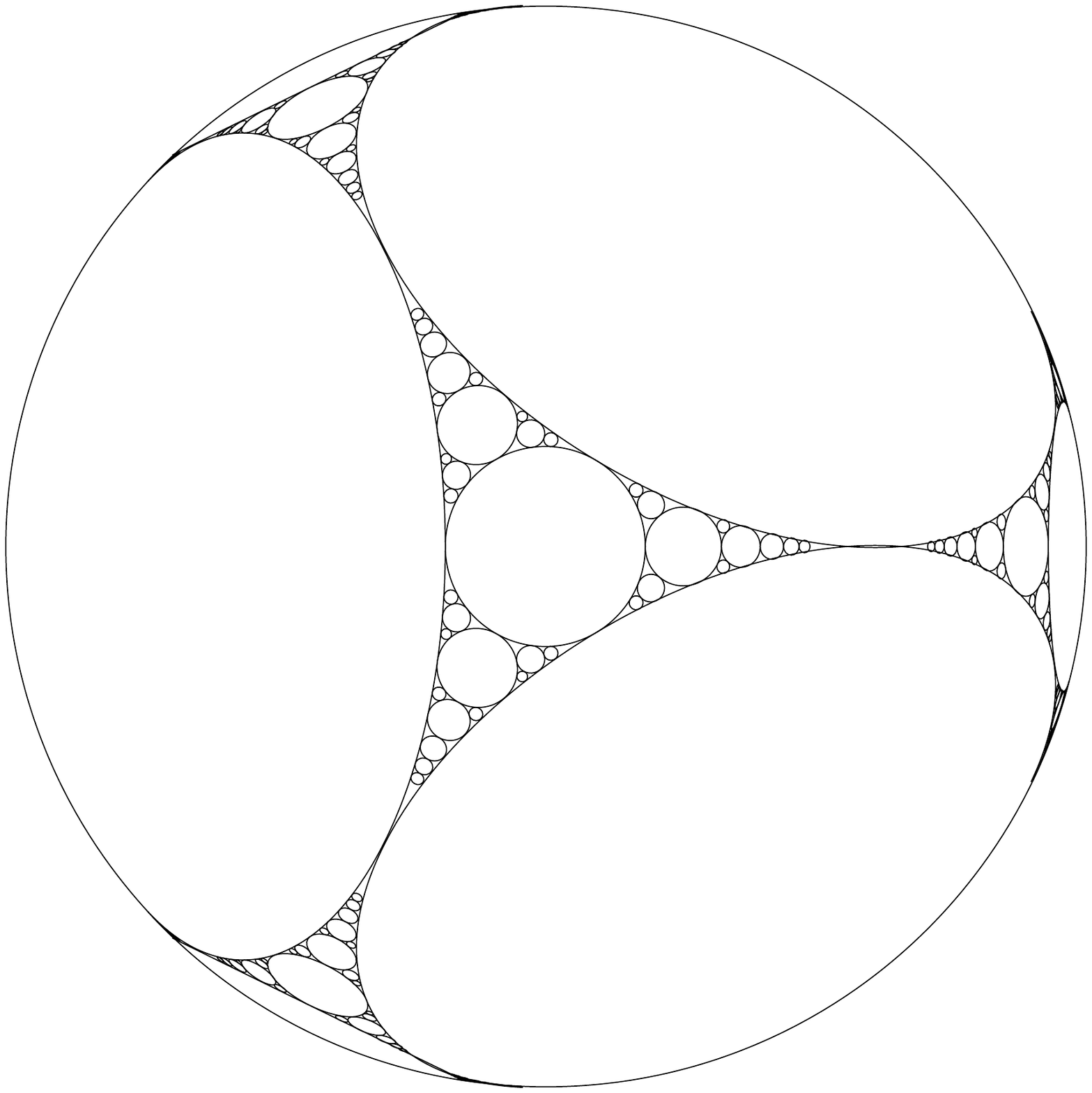}
    \includegraphics[width=5cm]{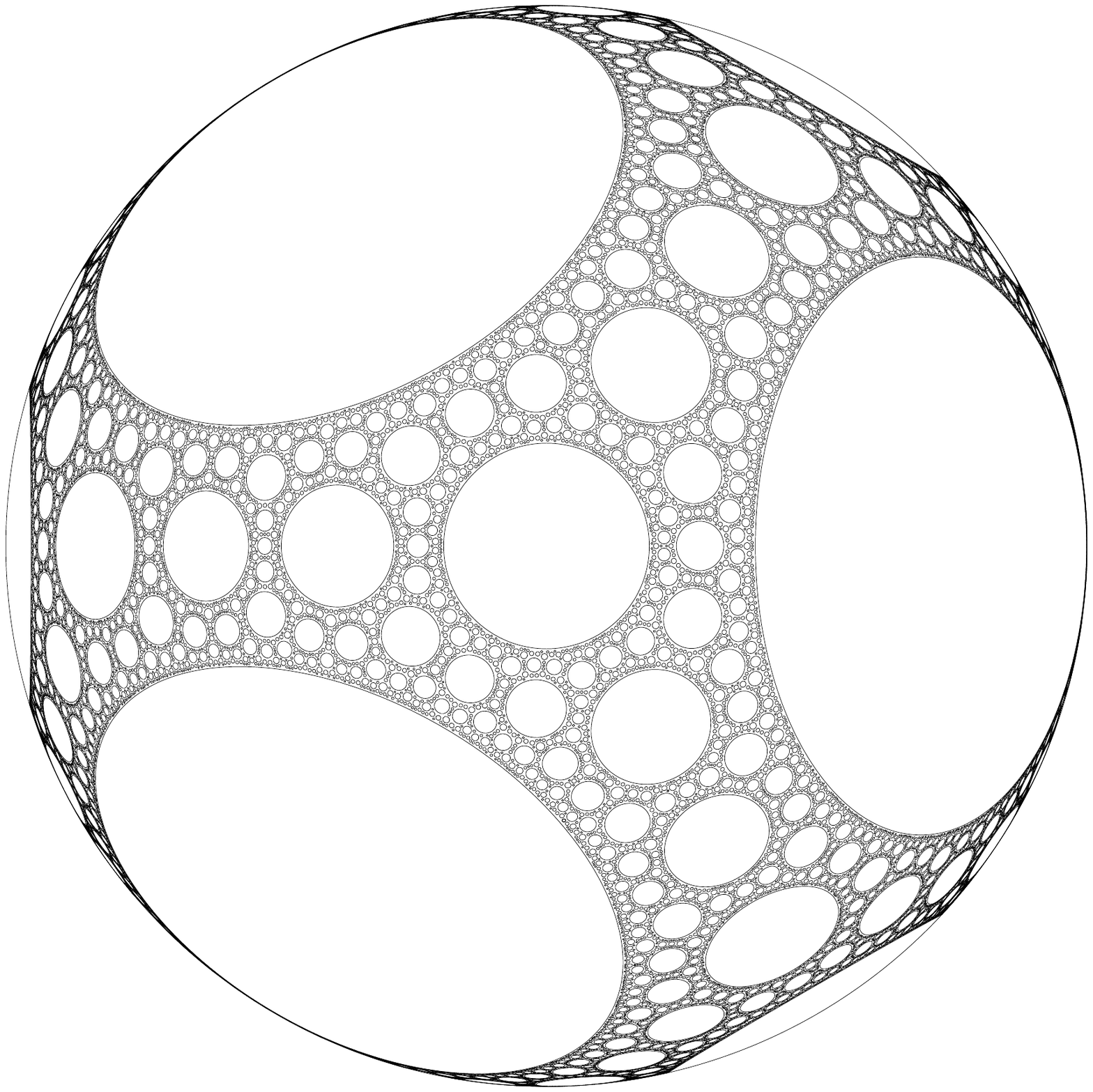}
     \caption{Apollonian packing and Sierpinski curve on the sphere (by C. McMullen)}
    \label{ASsphere}
 \end{center}
\end{figure}

Let $G$ denote the orientation preserving isometry group of $\mathbb B$.
\begin{Thm}[\cite{OhShahsphere}]
Let $\G$ be a non-elementary
geometrically finite discrete subgroup of $G$ and
 $\mathcal P=\cup_{i\in I} \Gamma(C_i)$ be an infinite,
  locally finite, and  $\G$-invariant circle packing on the sphere ${\mathbb S^2}$ with finitely many $\G$-orbits.

Suppose one of the following conditions hold:
\begin{enumerate}
\item $\G$ is convex co-compact;
\item  all circles in $\P$ are mutually disjoint;
\item $\cup_{i\in I} C_i^\circ \subset \Omega(\G)$
where $C_i^\circ$ denotes the interior of $C_i$.
% and $\Omega(\G):=\hat \c-\Lambda(\G )$ is the domain of discontinuity for $\G$.
\end{enumerate}
Then for any Borel subset $E\subset \mathbb S^2$ whose boundary is of zero Patterson-Sullivan measure,
 $$N_T(\P, E)
\sim
 \frac{ \sk_\G(\P)\cdot \nu_o (E) }{\delta_\G \cdot |m^{\BMS}_\G|} \cdot
  (2T)^{\delta_\G}  \quad\text{as $T\to \infty$}$$
 where $0<\sk_\G(\P)<\infty$ is defined in Def. \ref{sk}.
\end{Thm}

\section{Integral Apollonian packings: Primes and Twin primes}
A circle packing $\P$ is called {\it integral} if the curvatures of
all circles in $\P$ are integral.
One of the special features of Apollonian circle packings is the
abundant existence of {\it integral Apollonian circle packings}.

Descartes noted in 1643 (see \cite{Coxeter1968}) that
a quadruple $(a,b,c,d)$ of real numbers can be realized as curvatures of four mutually tangent circles
in the plane
(oriented so that their interiors are disjoint)
if and only if it satisfies
\begin{equation}\label{Des}
2(a^2+b^2+c^2+d^2)-(a+b+c+d)^2 =0.
%Q(a,b,c,d)=0
\end{equation}

% where
% \begin{equation}\label{des}
% Q(a,b,c,d)=2(a^2+b^2+c^2+d^2)-(a+b+c+d)^2 \end{equation} is the so-called Descartes quadratic form.

\begin{figure}
\begin{center}
\includegraphics[height=4cm] {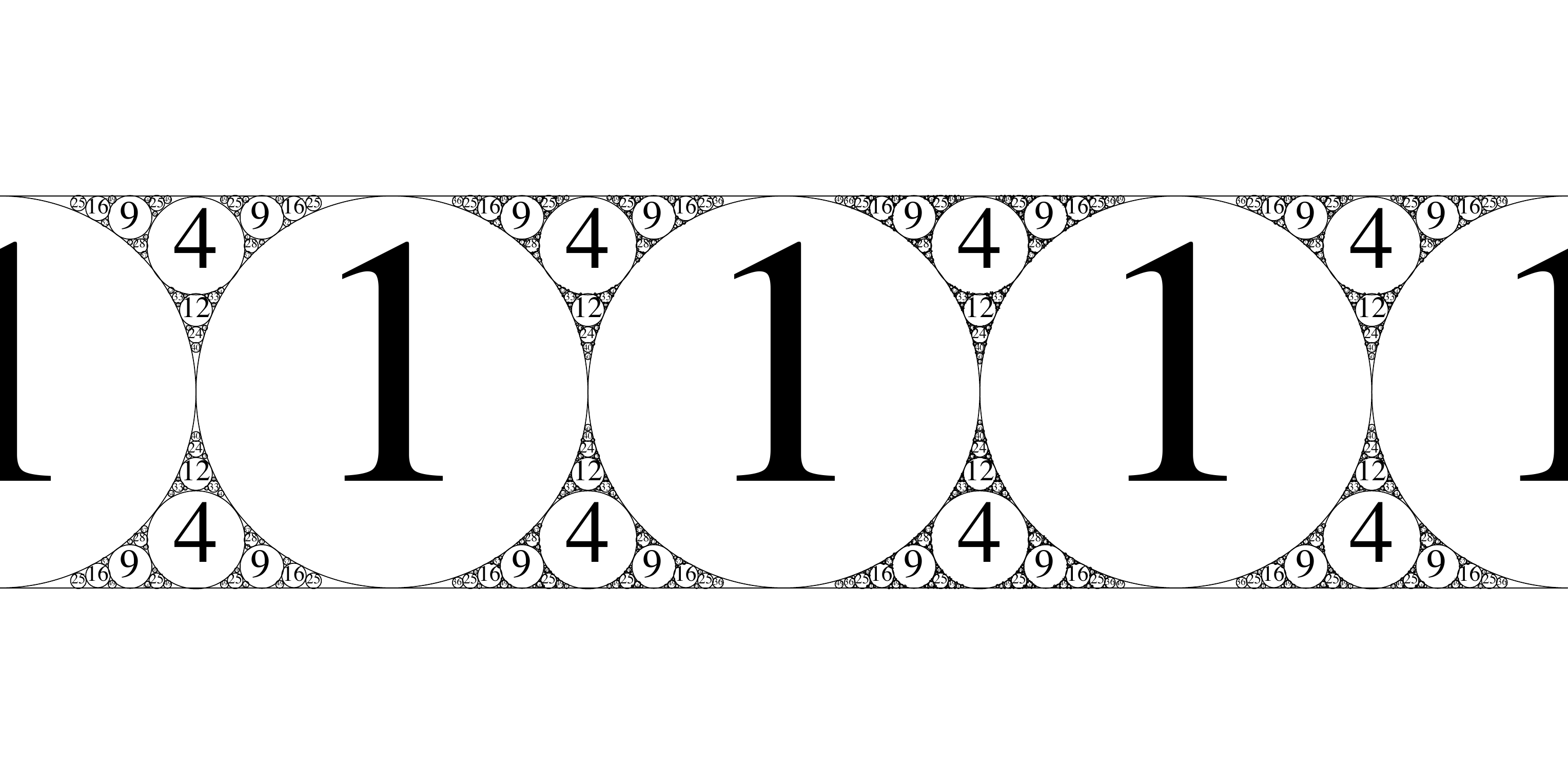}
\end{center}
\caption{An Apollonian circle packing between two parallel lines.}
\label{fpar}
\end{figure}
Usually referred to as the Descartes circle theorem,
this theorem implies that if the initial four circles in an Apollonian circle packing $\P$
in the plane have integral curvatures,
then $\P$ is an integral packing, as observed by
Soddy in 1937 \cite{Soddy1937}. The Descartes circle theorem provides an integral Apollonian packing
for every integral solution of the quadratic equation \eqref{Des} and indeed
there are infinitely many distinct integral Apollonian circle packings.

Let $\P$ be an integral Apollonian circle packing. We can deduce from the existence of the
lower bound
for the non-zero curvatures in $\P$ that
such $\P$ is either bounded or lies between two parallel lines.
We assume that $\P$ is primitive, that is,
the greatest common divisor of curvatures is one.

Calling a circle with a prime curvature a prime and
a pair of tangent prime circles a twin prime, Sarnak showed:
\begin{Thm}[\cite{SarnakToLagarias}]\label{sar}
There are infinitely many primes, as well as twin primes, in $\P$.
\end{Thm}

%%\begin{Rmk}{\rm
%%It would be interesting to know the infinitude of primes (and of twin primes)
%%in each triangular region of $\P$.
%%}\end{Rmk}

For $\mathcal P$ bounded,
denote by $\pi^{\mathcal P}(T)$ the
number of  prime circles in $\mathcal P$  of curvature at
most $T$, and  by $\pi^{\mathcal P}_2(T)$
 the number of twin prime circles in $\mathcal P$ of
 curvatures at most $T$.
 For $\mathcal P$ congruent to the packing in Fig. \ref{fpar},
  we alter the definition of $\pi^\P (T)$ and $\pi^{\mathcal P}_2(T)$ to count prime circles in a fixed period.
Sarnak showed \cite{SarnakToLagarias} that
  $$\pi^{\mathcal P}(T)\gg\frac{T}{{(\log T)}^{3/2}} .$$

%Let $\Gamma$ be a finitely generated subgroup of $\SL_2(\z)$
Recently Bourgain, Gamburd and Sarnak (\cite{BourgainGamburdSarnak2008}
and \cite{BourgainGamburdSarnak2008-p}) obtained a uniform spectral gap
for the family of
 congruence subgroups $\Gamma(q)=\{\gamma\in \G: \gamma\equiv 1 \;\; (\text{mod } q)\}$, $q$ square-free,
 of any finitely generated subgroup $\G$ of $\SL_2(\z)$ provided $\delta_{\G}>1/2$.
This theorem extends to a Zariski dense subgroup $\G$ of  $\SL_2(\z [i])$ and its congruence subgroups
over square free ideals of $\z[i]$ if $\delta_\G>1$.

Denoting by $Q$ the Descartes quadratic form
$$Q(a,b,c,d)=2(a^2+b^2+c^2+d^2)-(a+b+c+d)^2 ,$$
the approach in \cite{KontorovichOh} for counting circles in Apollonian circle
packings which are either bounded or between two parallel lines is
based on the interpretation of such circle counting problem into
the counting problem for $w\G\cap B_T^{\max}$
where $\G<O_Q(\z)$ is the so-called
Apollonian group,  $w\in \z^4$ with $Q(w)=0$ and
$B_T^{\max}$ denotes the maximum norm ball in $\br^4$.

Using the spin double cover $\op{Spin}_Q \to \SO_Q$ and the isomorphism $\op{Spin}_Q(\br)=\SL_2(\c)$,
we use the aforementioned result of Bourgain, Gamburd and Sarnak to obtain
a smoothed counting for $w\Gamma(q) \cap B_T$ with a uniform error term
for the family of square-free congruence subgroups $\Gamma(q)$'s where
$B_T$ is the Euclidean norm ball.
 This is a crucial ingredient for the Selberg's upper bound sieve, which is used to prove
 the following:

\begin{Thm}[\cite{KontorovichOh}]\label{prime} As $T\to \infty$,
$$ \pi^{\mathcal P}(T)\ll  \frac{T^{\alpha}}{\log T},\quad\text{and}\quad
 \pi^{\mathcal P}_2(T)\ll  \frac{T^{\alpha}}{(\log T)^2} $$
 where $\alpha=1.30568(8)$ is the residual dimension of $\P$.
 \end{Thm}

\begin{Rmk}{\rm
\begin{enumerate}
\item
Modulo $16$, the Descartes equation \eqref{Des} has no solutions unless two of the curvatures
are even and the other two odd. In particular,
there are no ``triplet primes'' of three mutually tangent circles, all having odd prime curvatures.

\item
We can also use the methods in
\cite{KontorovichOh} to give lower bounds for  almost primes in a packing. A circle in $\P$ is called {\it $R$-almost prime} if its curvature is the product of at most $R$ primes. Similarly, a pair of tangent circles is called $R${\it-almost twin prime} if both circles are $R$-almost prime.
Employing Brun's combinatorial sieve, our methods show the existence of $R_1, R_2>0$ (unspecified) such that the number of $R_1$-almost prime
circles in $\mathcal P$ whose curvature is at most $T$
is $\asymp \frac{T^{\alpha}}{\log T},$\footnote{
By $f(T)\asymp g(T)$, we  mean
$g(T)\ll f(T)\ll g(T)$.}
 and that the number of
pairs of  $R_2$-almost twin prime circles whose curvatures are at most $T$
is $\asymp \frac{T^{\alpha}}{(\log T)^2}$.
%\footnote{
%The infinitude and in fact Zariski density of $R$-almost prime and twin prime circles is obtained in
%\cite{BourgainGamburdSarnak2008}. The difference between our statements and theirs is that we can count with an archimedean norm. Neither method can specify $R_1$ or $R_2$, due to the lack of an explicit spectral gap, see \cite{BourgainGamburd2007}.
%}

\item A suitably modified version of Conjecture 1.4 in \cite{BourgainGamburdSarnak2008},
 a generalization of Schinzel's hypothesis,
implies that for some $c, c_2>0$,
$$\pi^{\mathcal P}(T) \sim c \cdot \frac{T^{\alpha}}{\log T}\quad\text{and}\quad
\pi^{\mathcal P}_2(T) \sim c_2 \cdot \frac{T^{\alpha}}{(\log T)^2}.
$$
%(see the discussion in \cite[Ex D]{BourgainGamburdSarnak2008}).
The constants $c$ and $c_2$ are detailed in \cite{FuchsSanden}.
\item Recently Bourgain and Fuchs \cite{BourgainFuchs2010} showed that in a given bounded integral Apollonian
packing $\P$, the growth of the number
of {\it distinct} curvatures at most $T$ is at least $c\cdot T$ for some $c>0$.
\item The spherical Soddy-Gossett theorem says (see \cite{LagariasMallowsWilks}) that
 the quadruple $(a,b,c,d)$ of spherical curvatures of four mutually tangent circles
in $\P$ satisfies
$$2(a^2+b^2+c^2+d^2)-(a+b+c+d)^2=-4 .$$
This theorem implies again that there are infinitely many {\it integral} spherical Apollonian circle packings, that is, the spherical curvature of every circle is integral.
It will be interesting to have results analogous to
Theorems  \ref{sar} and \ref{prime} for integral spherical Apollonian packings.
\end{enumerate}} \end{Rmk}

\section{Equidistribution in geometrically finite hyperbolic manifolds}\label{br}

%We now explain our equidistribution results
 %which are the main ingredients in proving counting theorems in the previous section.
  Let $G$ be the identity component of the group of isometries of
  $\bH^n$ and $\G<G$ be a non-elementary geometrically finite discrete
   subgroup.

We have discussed that the Bowen-Margulis-Sullivan measure is
a finite measure on the unit tangent bundle $\T^1(\G\ba \bH^n)$ which
is mixing for the geodesic flow.
Another measure playing an important role in studying the dynamics
of flows on $\T^1(\G\ba \bH^n)$ is the following Burger-Roblin measure.

\bigskip

\noindent{\bf {Burger-Roblin measure}:}
 The Burger-Roblin measure $m^{\BR}_\G$
is the induced measure on $\T^1(\G\ba \bH^n)$ of the following $\G$-invariant measure on $\T^1(\bH^n)$:
 $$d \tilde m^{\BR}(u)=e^{(n-1) \beta_{u^+}(x,
\pi(u))}\;
 e^{\delta_{\G} \beta_{u^-}(x, \pi(u)) }\; dm_x(u^+) d\nu_x(u^-) dt $$
where $m_x$ denotes the probability measure
 on the boundary
$\partial_\infty(\bH^n)$ invariant under  the maximal compact subgroup $\op{Stab}_G(x_0)$.
 %a $G$-invariant conformal density on the boundary
%$\partial_\infty(\bH^n)$ of dimension $(n-1)$,
%unique up to homothety.
%Explicitly, fix any $x_0\in \bH^n$ and let $m_{x_0}$ denote the unique probability measure on
%$\partial_\infty(\bH^n)$ invariant under the maximal compact subgroup $\op{Stab}_G(x_0)$.
For any $x$ and $ x_0\in \bH^n$,
we have $dm_x(\xi):=e^{-(n-1)\beta_{\xi} (x,x_0)}dm_{x_0}(\xi)$ and
it follows that
 this definition of $m^{\BR}_\G$ is independent of the choice of $x\in \bH^n$.
%We denote by $m^{\BR}$ the induced measure on $\T^1(X)$.}
%\end{Def}
%In particular, each $m_x$ is a finite measure on
%$\partial_\infty(\bH^n)$ invariant under the stabilizer of $x$ in $G$,
%and the family $\{m_x\}$ satisfies that
% for any $x,y\in \bH^n$, $\xi\in \partial_\infty(\bH^n)$ and $g\in
%G$,
%$$g_*\nu_x=\nu_{g x}\quad\text{ and}\quad \frac{d\nu_y}{d\nu_x}(\xi)=e^{-(n-1)\beta_{\xi} (y,x)}. $$

Burger \cite{Burger1990} showed that for a convex cocompact hyperbolic surface with
$\delta_\G$ at least $1/2$,
 this is a unique ergodic horocycle invariant measure up to homothety.
  Roblin \cite{Roblin2003} extended
 Burger's result in much greater generality, for instance,
 including all non-elementary geometrically finite hyperbolic manifolds.

The name of the Burger-Roblin measure
  was first suggested by Shah and the author in \cite{KontorovichOh} and
 \cite{OhShahGFH} in recognition of this
  important classification result.
  % is not directly
 %used in our arguments below, but it certainly
 %explains a reason behind our result that the asymptotic distribution of expanding horospheres
 %is
 %described by the Burger-Roblin measure.

We note that the total mass $|m^{\BR}_{\G}|$ is finite only when $\delta_\G=n-1$ (or
equivalently only when $\G$ is a lattice in $G$) and is supported
on the set $\{u\in \op{T}^1(\G\ba \bH^n): u^-\in \LG\} .$

\bigskip

    Let $S^\dag \subset \T^1(\bH^n)$ be one of the following:
   \begin{enumerate} \item an unstable horosphere;
\item the oriented unit normal bundle of a codimension
one totally geodesic subspace of $\bH^n$
     \item the set of outward
 normal vectors to a (hyperbolic) sphere in $\bH^n$.
\end{enumerate}

% Assuming that
% the image $p(\tS)$ is closed in the unit tangent bundle
% $\T^1(\G\ba \bH^n)$,
 %the main ergodic ingredient needed in proving Theorem \ref{m1} is
 %to understand the asymptotic
 %distribution of $p(\tS)$ when pushed by the geodesic flow $g^r$ as $r\to \infty$.

%To be precise, we recall that generalizing the work of Patterson
%for $n=2$ \cite{Patterson1976}, Sullivan \cite{Sullivan1979}
%constructed a $\G$-invariant conformal density $\{\nu_x:x\in
%\bH^n\}$ of dimension $\delta_\G$ on $\Lambda(\G)$.
We consider the following measures on $\op{Stab}_\G( S^\dag)\ba S^\dag:$
$$d\mu^{\Leb}_{S^\dag}(s)=e^{(n-1)\beta_{s^+}(x,\pi(s))}dm_x(s^+),\quad
d\mu^{\PS}_{S^\dag}(s)=e^{\delta_\G\beta_{s^+}(x,\pi(s))}d\nu_x(s^+)$$
for any $x\in \bH^n$.

Denote by $p$ the
 canonical projection $\T^1(\bH^n)\to \T^1(\G\ba \bH^n)=\G\ba
 \T^1(\bH^n)$.
\begin{Thm}[\cite{OhShahGFH}] \label{main}
 For $\psi\in C_c(\T^1(\G\ba \bH^n))$ and any relatively compact subset $\mathcal O\subset
 p(S^\dag)$ with $\mu^{\PS}_{S^\dagger}(\partial(\mathcal O))=0$,
\begin{equation*} e^{(n-1-\delta_\G)t}\cdot
\int_{\mathcal O}\psi(g^t(s))\; d\mu^{\Leb}_{S^\dag}(s)
\sim \frac{\mu_{S^\dag}^{\PS}(\mathcal O_*)}{\delta_\G\cdot |m^{\BMS}_\G|}
 \cdot
m^{\BR}_\G(\psi)\quad\text{as $t\to \infty$}
\end{equation*}
where $$\mathcal O_*=\{s\in\mathcal O: s^+\in \Lambda(\G)\}.$$
\end{Thm}

\begin{Def} {\rm For a hyperbolic subspace $S=\bH^{n-1}\subset \bH^n$,
we say that a parabolic fixed point $\xi \in \Lambda(\G)\cap \partial_\infty(\bH^{n-1})$ of $\G$
is {\it internal} if any parabolic element $\gamma\in \G$
fixing $\xi$ preserves $\bH^{n-1}$.}
\end{Def}

Recalling the notation $\pi$ for the canonical projection from $\T^1(\bH^n)$ to $\bH^n$,
 we set $S=\pi(S^\dag)$.
\begin{Thm}[\cite{OhShahGFH}] \label{mainer} We assume that the projection map
 $\op{Stab}_\G(S)\ba  S
\to \G \ba \bH^n$
is proper.
 In the case when ${S}$ is a codimension one totally geodesic
subspace, we also assume that every parabolic fixed point of $\G$ in the boundary
 of $S$ is internal.

 For $\psi\in C_c(\T^1(\G\ba \bH^n))$,
\begin{equation*} e^{(n-1-\delta_\G)t}\cdot
\int_{p({S^\dag})}\psi(g^t(s))\; d\mu^{\Leb}_{{S^\dag}}(s)
\sim \frac{\mu_{{S^\dag}}^{\PS}(S^\dag_*)}{\delta_\G\cdot |m^{\BMS}_\G|}
 \cdot
m^{\BR}_\G(\psi)\quad\text{as $t\to \infty$}
\end{equation*}
where $$S^\dag_*=\{s\in p({S^\dag}): s^+\in \Lambda(\G)\}.$$

We have $0\le \mu_{{S^\dag}}^{\PS}(S^\dag_*)<\infty,$
and $\mu_{{S^\dag}}^{\PS}(S^\dag_*)=0$ may happen only when $S$ is totally geodesic.
\end{Thm}

It can be shown by combining results of
\cite{Dani1979} and \cite{MargulisICM1991} that
in a finite volume space $\G\ba \bH^n$,
the properness of the projection map
 $\op{Stab}_\G  (S) \ba  S
\to \G \ba \bH^3$
implies that $\op{Stab}_\G( S)\ba S$ is of finite volume as well, except
for the case when $n=2$ and $S$ is a proper geodesic in $\bH^2$ connecting two
parabolic fixed points of a lattice $\G <\PSL_2(\br)$.

When both $\G\ba \bH^n$ and $\op{Stab}_\G( S)\ba S$ are of finite volume,
we have $n-1=\delta_\G$ and both $m^{\BMS}_\G$ and $m^{\BR}_\G$
are finite invariant measures and $\mu^{\PS}_{S^\dag}= \mu^{\Leb}_{S^\dag}$ (up to a constant multiple).
In this case, Theorem \ref{mainer} is
 due to Sarnak \cite{Sarnak1981} for the closed horocycles for $n=2$.
 The general case is due to
  Duke, Rudnick and Sarnak \cite{DukeRudnickSarnak1993}
and
Eskin and McMullen \cite{EskinMcMullen1993} gave a simpler proof of Theorem \ref{mainer},
based on the mixing property of the geodesic flow of a finite volume
hyperbolic manifold. The latter proof, combined with a strengthened version
of the wavefront lemma \cite{GorodnikShahOhstrongwave}, also works for proving
Theorem \ref{main}.
We remark that the idea of using mixing in this type of problem
  goes back to the 1970 thesis of Margulis \cite{Margulisthesis}
(see also \cite[Appendix]{KleinbockMargulis1996}).
Eskin, Mozes and Shah \cite{EskinMozesShah1996} and Shah \cite{Shah1996}
provided yet another different proofs using the theory of unipotent flows.
When both $\G\ba \bH^n$ and $\op{Stab}_\G( S)\ba S$ are of finite volume,
Theorem \ref{main} easily implies Theorem \ref{mainer} but not conversely.

In the case when $S^\dag$ is a horosphere, Theorem \ref{main} was obtained
 in \cite{Roblin2003}, and Theorem \ref{mainer} was proved
in \cite{KontorovichOh} when $\delta_\G >(n-1)/2$ with a different
interpretation of the main term.

\begin{Rmk}
{\rm
\begin{enumerate}
\item The condition on the internality of all parabolic fixed points of $\G$ in the boundary
of $S$ is crucial, as $\mu_{S^\dagger}^{\PS}(S^\dag_*)=\infty$ otherwise.
This can already be seen in the level of a lattice: take $\G=\SL_2(\z)$
 and let $S$ be the geodesic connecting $0$ and $\infty$ in the upper half space model.
 Then any upper triangular matrix in $\G$ fixes $\infty$ but does not stabilize $S$.
Indeed the length of the image of $S$ in $\G\ba \bH^2$ is infinite.

\item In proving Theorem \ref{m1}, we count circles in $\c$ by counting the
 corresponding Euclidean hemispheres in $\bH^3$. As the Euclidean hemispheres
 are totally geodesic hyperbolic planes, this amounts to understanding the distribution
 of a $\Gamma$-orbit of a totally geodesic hyperbolic plane in $\bH^3$.
 The equidistribution theorem we use here is Theorem \ref{mainer} for $S$ a hyperbolic plane.

\item More classical applications of the equidistribution theorem
such as Theorem \ref{mainer} can be found in the point counting problems of
 $\G$-orbits in various spaces.

For a $\G$-orbit in the hyperbolic space $\bH^n$, the orbital counting
in Riemannian balls was obtained
 Lax-Phillips \cite{LaxPhillips} for $\delta_\G>\frac{n-1}{2}$ and by
     Roblin \cite{Roblin2003} in general.

Extending the work of Duke, Rudnick and Sarnak \cite{DukeRudnickSarnak1993} and
of Eskin and McMullen \cite{EskinMcMullen1993} for $\G$ lattices,
 we obtain in \cite{OhShahGFH}, for any geometrically finite group $\G$ of $G$,
the asymptotic of the number of vectors of norm at most $T$ lying in a discrete orbit $w \Gamma$
of a quadric $$F(x_1, \cdots, x_{n+1})=y$$
for a real quadratic form $F$ of signature $(n,1)$ and any $y\in \br$
(when $y>0$, there is an extra assumption on $w$ not being $\G$ strongly parabolic.
See \cite{OhShahGFH} for details).
When $y=0$ and $n=2, 3$, special cases of this result were obtained in
\cite{Kontorovich2007}, \cite{KontorovichOh} and \cite{KontorovichOh2008}
under the condition $\delta_\G>(n-1)/2$.
Based on the Descartes circle theorem, this result in
\cite{KontorovichOh} was used to prove Theorem \ref{apo} for the bounded Apollonian packings.
In \cite{LimOhboundary}, a $\G$-orbit in the geometric boundary
 is shown to be equidistributed with respect to the Patterson-Sullivan measure, extending
 the work \cite{GorodnikOh2007} for the lattice case.
%Applications to Pythagorean triples and to Apollonian
%circle packings are described in there.

\item For $\psi\in C_c(\G\ba \bH^n)$,
we have
$$m^{\BR}_\G(\psi)=\la \psi, \phi_\G\ra :=\int_{\G\ba \bH^n} \psi(x) \cdot \phi_\G (x)\;  dm^{\Leb}(x)$$
where $\phi_\G(x)=|\nu_x|$ is the positive eigenfunction of the Laplace
operator on $\G\ba \bH^n$ with eigenvalue $\delta_\G(n-1-\delta_\G)$
and $$dm^{\Leb}(u)=e^{(n-1) \beta_{u^+}(x,
\pi(u))}\;
 e^{(n-1) \beta_{u^-}(x, \pi(u)) }\; dm_x(u^+) dm_x(u^-) dt $$
for any $x\in \bH^n$.
Hence Theorem \ref{mainer} says that for $\psi\in C_c(\G\ba \bH^n)$,
\begin{equation}\label{er} e^{(n-1-\delta_\G)t}\cdot
\int_{p({S^\dag})}\psi(\pi(g^t(s)))\; d\mu^{\Leb}_{{S^\dag}}(s)
\sim \frac{\mu_{{S^\dag}}^{\PS}(S^\dag_*)}{\delta_\G\cdot |m^{\BMS}_\G|}
 \cdot
\la \psi, \phi_\G\ra \quad\text{as $t\to \infty$} .\end{equation}

When $\delta_\G >(n-1)/2$, $\phi_\G\in L^2(\G\ba \bH^n)$ and its eigenvalue
$\delta_\G(n-1-\delta_\G)$
 is isolated in the $L^2$-spectrum of the Laplace operator \cite{LaxPhillips}.
  It will be desirable to obtain a rate of convergence
  in \eqref{er} in terms of the spectral gap of $\G$ in such cases. For $\G$ lattices, it was achieved in
  \cite{DukeRudnickSarnak1993} for $p(S)$ compact and in \cite{BenoistOh2008}
in general.  This was done in the case of a horosphere
 in \cite{KontorovichOh}, which was the main ingredient in the proof
 of Theorem \ref{prime}. It may be possible to extend the methods of \cite{KontorovichOh} to obtain an error term in general.
\end{enumerate}
}
\end{Rmk}

% The
%submanifold $p(\tS)$ is not necessarily compact or not even of finite volume in general.
%However we show based on the geometry at infinity of a
%geometrically finite hyperbolic manifold that the integral on the
%left hand side of \eqref{mainer} takes place only over a compact
%subset and that $\tS_*^+$ is always compact. It is in this last
%step that we use the assumption that $\pi(\tS)$ does not meet any
%parabolic fixed point for $\G$ in the boundary when $\tS$ is based
%on a totally geodesic hypersurface.

%As we do not impose any condition on $\delta_\G$ but
%to be positive, such a square-integrable base eigenfunction $\phi_\G$ does not exist any more
%and hence we cannot rely on the spectral argument on $L^2(\G\ba
%\bH^n)$. Instead we use
%the dynamics of the geodesic flow on the unit tangent bundle $\T^1(\G\ba \bH^n)$.

\section{Further remarks and questions}
Let $G$ be the identity component of the group
of isometries of $\bH^n$ and
 $\G$ be a geometrically finite group.
We further assume that $\G$ is Zariski dense in $G$ for discussions in this section.
 When we identify $\bH^n$ with $G/K$
 for a maximal compact subgroup $K$, the unit tangent bundle $\T^1(\bH^n)$
 can be identified with $ G/M$ where $M$ is the centralizer
in $K$ of a Cartan subgroup, say, $A$, whose multiplication on the right corresponds to
the geodesic flow.
The frame bundle of $\bH^n$ can be identified with $G$ and
the frame flow on the frame bundle is given by the multiplications
by elements of $A$ on the right.

We have stated the equidistribution results in section \ref{br}
 in the level of the unit tangent bundle
$\T^1(\G\ba \bH^n)$.
As the frame bundle is a homogeneous space of $G$ unlike the unit tangent bundle,
it is much more convenient to work in the frame bundle.
Fortunately, as observed in \cite{FlaminioSpatzier}, the frame flow is mixing on
$\G\ba G$ with respect to the lift from $\G\ba G/M$ to $\G\ba
G$ of the Bowen-Margulis-Sullivan measure.
Using this, we can extend Theorems \ref{main} and \ref{mainer} to the level of the frame
bundle $\G\ba G$.
It seems that the classification theorem of Burger and Roblin
 can also be extended: for a horospherical group $N$,
any locally finite $N$-invariant ergodic measure on $\G\ba G$ is either supported on a closed
$N$-orbit or the lift of the Burger-Roblin measure (we caution here that
 a locally finite $N$-invariant measure
supported on a closed $N$-orbit need not be a finite measure unlike the $\G$-lattice cases).

In analogy with Ratner's theorem \cite{Ratner1991}, \cite{Ratner1991t}, we
propose the following problems:
let $U$ be a one-parameter unipotent subgroup, or more generally a subgroup generated by unipotent
one parameter subgroups of $G$:
\begin{enumerate}
\item{[Measure rigidity]} Classify all locally finite Borel
$U$-invariant ergodic measures on $\Gamma \ba G$.
\item{[Topological rigidity]} Classify the closures of
$U$-orbits in $\Gamma \ba G$.
\end{enumerate}

We remark that as $G=\op{SO}(n,1)$ (up to a local isomorphism) in our set-up,
the above topological rigidity for $\G$ lattices was also obtained by Shah (\cite{Shah1991},
\cite{Shah1991Trieste}) based on the apporoach of Margulis (\cite{Margulis1986}, \cite{Margulis1989})
and of Dani and Margulis \cite{DaniMargulis1990}.

Both questions are known for $n=2$ due to Burger \cite{Burger1990} and
Roblin \cite{Roblin2003}, as in this case, there is only one unipotent one-parameter
subgroup up to conjugation, which gives the horocycle flow.
Shapira used them to prove equidistribution for non-closed horocycles \cite{Schapira2005}.

It may be a good idea to start with a sampling case when $G=\SL_2(\c)$, $U=\SL_2(\br)$ and
 $\G <G$ Zariski dense and geometrically finite.
\begin{enumerate}
\item Are there any locally finite $\SL_2(\br)$-invariant ergodic measure on $\G\ba \SL_2(\c)$
besides the Haar measure (=the $\SL_2(\c)$-invariant measures)
 and the $\SL_2(\br)$-invariant measures supported on closed $\SL_2(\br)$ orbits?
 \item  Is every non-closed $\SL_2(\br)$-orbit dense in $\G\ba \SL_2(\c)$ ?
\end{enumerate}

It seems that the answers are {\it no} for (1) and {\it yes} for (2).

\bibliographystyle{plain}
%\bibliography{Ohbibliog}
%\end{document}

\end{document}